\documentclass[12pt]{article}
\usepackage{mathrsfs}
\usepackage[pdftex]{graphicx}
\usepackage{amsmath,amssymb}
\usepackage{amsfonts}
\usepackage{txfonts}

\newtheorem{theorem}{Theorem}
\newtheorem{lemma}{Lemma}

\newtheorem{definition}{Definition}[section]
\newtheorem{corollary}{Corollary}[section]

\newcommand\figcaption{\def\@captype{figure}\caption}

\begin{document}
\title{A KAM-Theorem for Persistence of Quasi-periodic Invariant Tori in Bifurcation Theory of Equilibrium Points
\thanks{This work is supported  by the
NNSF(11371132,11671392) of China, by Key Laboratory of High Performance Computing and Stochastic Information Processing. $^a$email:
lixuemei$\_1$@sina.com, $^b$ email: zaijiu@amss.ac.cn} }
\author{ Xuemei Li$^a$ and Zaijiu Shang$^b$ \\
\small $^{a}$ Key Laboratory of High Performance Computing and  Stochastic Information Processing,\\
 \small Department of Mathematics, Hunan Normal University, Changsha, Hunan 410081, P. R. China\\
\small $^{b}$ 1. HLM, Academy of Mathematics and Systems Science, Chinese Academy of Sciences, Beijing 100190, China\\
\small 2. School of Mathematical Sciences, University of Chinese Academy of Sciences, Beijing 100049, China
\\}
\date{}
\maketitle \vskip 0.3in

{\bf Abstract.} In this paper, we establish a KAM-theorem for ordinary differential equations with finitely differentiable vector fields and multiple degeneracies. The theorem can be used to deal with the persistence of quasi-periodic invariant tori in multiple Hopf and zero-multiple Hopf bifurcations, as well as their subordinate bifurcations, of equilibrium points of continuous dynamical systems.

\vskip 0.2in

\renewcommand{\theequation}{\thesection.\arabic{equation}}
\section*{ 1. Introduction}
\setcounter{section}{1}\setcounter{equation}{0}

 To study the bifurcations of equilibria of a system of differential equations (ODEs, PDEs and functional differential equations), one usually reduces such a system to a lower-dimensional one on the center manifold by the Center Manifold Theorem.  Possibly the reduced system  is finitely differentiable even if the original system is analytic. When the equilibrium is partially elliptic and the normal form of the reduced subsystem has a normal form of Birkhoff type on the center manifold, then the truncated normal form may possess quasi-periodic invariant tori. In this case, a question arises naturally: does the original system (equivalently, the reduced system on the center manifold) have quasi-periodic invariant tori with the same dimension ? This problem can be discussed by KAM theory and a careful study leads us to consider the existence of quasi-periodic tori of the following system
\begin{equation}\label{1.1}
\left\{\begin{array}{rl}
\dot I_1 & = \varepsilon^{q_1} [A_1(\xi,\varepsilon)I_1+\varepsilon^{q_2}g_1(I,\varphi; \xi,\varepsilon)]\\
\dot I_2 & = \varepsilon^{q_3} [A_2(\xi,\varepsilon)I_2+\varepsilon^{q_4}g_2(I,\varphi; \xi,\varepsilon)]\\
\dot{\varphi}_1 & = \varepsilon^{q_5}[\omega_1(\xi,\varepsilon)+\varepsilon^{q_6}g_3(I,\varphi; \xi,\varepsilon)]\\
\dot{\varphi}_2 & =\omega_2(\xi,\varepsilon)+\varepsilon^{q_7}g_4(I,\varphi; \xi,\varepsilon),
\end{array}
\right.
 \end{equation}
 where $I={\rm col}(I_1,I_2)\in \Omega\subset \mathbb{R}^{n_{11}}\times \mathbb{R}^{n_{12}}=\mathbb{R}^{n_1}, \varphi={\rm col}(\varphi_1,\varphi_2)\in \mathbb{T}^{n_{21}}\times \mathbb{T}^{n_{22}}=\mathbb{T}^{n_2}, q_j\geq 0(j=1,\cdots,7), \xi\in \Pi\subset \mathbb{R}^{n_3}$ is the bifurcation parameter, $\varepsilon$ is a small perturbation parameter.

When $g_j=0, j=1,\cdots,4$,
 $I=0$ represents the quasi-periodic torus of the integrable part of \eqref{1.1} which corresponds to the invariant torus of truncated normal forms. The aim of the present paper is to examine the persistence of the quasi-periodic torus $I=0$ under small perturbations (i.e., $g_j\neq 0,j=1,\cdots,4$). We meet some difficulties: the perturbation terms $g_j,j=1,\cdots,4$ are only finitely differentiable, there exist small frequencies, small twist and higher-order degeneracy in \eqref{1.1} and the number of parameter variables is possibly less than the dimension of tori. We need to tackle these difficulties in constructing a new KAM theorem for \eqref{1.1}.

 In the context of finitely differentiable perturbations, the study on the persistence of quasi-periodic invariant tori has originated from the work of Moser \cite{Mos62} on area-preserving mappings of an annulus, which was extended to dissipative vector fields in \cite{BMS76} based on smoothing operator technique. Another important method, which can relax the requirement for regularity of perturbations,  is to approximate a differentiable function by real analytic ones \cite{ Mos66, Russ70, Zehn75, Pos82, Sh00, CQ04, Alb07, Wag10}. R\"{u}ssmann proved an optimal estimate result on approximating a differentiable function by analytic ones. Following this approach Zehnder \cite{Zehn75} established a generalized implicit function theorem and applied it to the existence of parameterized invariant tori of nearly integrable Hamiltonian systems in finitely differentiable case, P\"{o}schel \cite{Pos82} showed that on a Cantor set, invariant tori of the perturbed Hamiltonian system form a differentiable family in the sense of Whitney. The results and ideas of Moser and P\"{o}schel are extended to the case of symplectic mappings by Shang \cite{Sh00} and to the case of lower dimensional elliptic tori by Chierchia and Qian \cite{CQ04}, respectively. Wagener \cite{Wag10} extended the modifying terms theorem of Moser \cite{Mos67} (i.e., introducing additional parameters) to finitely differentiable and Gevrey regular vector fields.

 The results mentioned above, except for \cite{CQ04}, were restricted to the case where the integrable part is analytic in coordinate variables as well as in parameters. The integrable part in \cite{CQ04} is assumed to be Lipschitz with respect to parameters and the frequency map to be a Lipschitz homeomorphism. Of course, if the unperturbed (integrable) part and the perturbation are both of class $C^l$ with $l>2n$ ($n$ is the number of degrees of freedom), it is reduced to the case where the integrable part is analytic and the perturbation is of class $C^l$ by regarding the initial values of action variables as parameters. The KAM theorems in \cite{BHS96, Wag10} can be applied to quasi-periodic bifurcations (i.e., bifurcations of quasi-periodic invariant tori). In this paper, we shall extend the result and method of P\"{o}schel \cite{Pos82} to the dissipative system \eqref{1.1} with degeneracies, and provide a convenient tool to investigate the persistence of quasi-periodic invariant tori in bifurcation theory of equilibrium points.

 The perturbation was assumed to be $C^{333}$ originally in the work of Moser \cite{Mos62} on area-preserving mappings of an annulus, and then was weakened to $C^5$ by R\"{u}ssmann \cite{Russ70} and  to $C^l (l>3)$ (meaning that the perturbation is of class $C^3$ and the derivatives of order 3 are H\u{o}lder continuous) by R\"{u}ssmann \cite{Russ83} and Herman \cite{Her83}, where a counterexample for $l<3$ was given. For improvements on weakening the regularity of perturbations in the Hamiltonian case we refer to \cite{Alb07} and references therein.

 The above mentioned results were proved under the so-called non-degeneracy conditions. In the context of degenerate KAM theory, i.e., if Kolmogorov's non-degeneracy or Arnold's iso-energetic non-degeneracy condition is violated, Arnol'd \cite{Arn63} established a properly degenerate KAM theorem (refined by \cite{Fej04, Chier10}) to deal with quasi-periodic motions in the planetary many body problem. In this case the integrable part does not depend on the full set of action variables, and the non-degeneracy conditions are imposed additionally on the averaged perturbation. The ideas of Arnol'd \cite{Arn63} were extended to the resonant torus case in \cite{CW99, LY03} and the normal zero-frequency case in \cite{GG05, HLY06, Gen07} for Hamiltonian systems and in \cite{BG01, LY12, Li16} for dissipative systems. Another method is to search for weaker non-degeneracy conditions concerning frequency maps, which have been studied in a series of papers, for example, by Bruno \cite{Bru92}, Cheng and Sun \cite{CS94}, R\"{u}ssmann \cite{Russ89, Russ01}, Han, Li and Yi \cite{HLY10} for finite dimensional Hamiltonian systems, and Bambusi, Berti and Magistrelli \cite{BBM11} for infinite dimensional case. The weaker non-degeneracy condition in \cite{CS94} is that the image of the frequency map in an open set includes a curved $C^{n+2}$ one-dimensional submanifold. R\"{u}ssmann \cite{Russ89, Russ01} pointed out that the weaker non-degeneracy condition means that the image of the frequency map does not lie in an $(n-1)$-dimensional linear subspace of $\mathbb{R}^n$ (this condition is also necessary in the analytic case).  An interesting and real analytic Hamiltonian of the form
$$H(x,y,\varepsilon)=h_0(y^{n_0})+\varepsilon^{m_1} h_1(y^{n_1})+\cdots + \varepsilon^{m_a}h_a(y^{n_a})+\varepsilon^{m_a+1}P(x,y,\varepsilon)$$
with the degeneracy involving several time scales was considered in  \cite{HLY10}. The degeneracy in \eqref{1.1} is somewhat similar to the one in \cite{HLY10}.

\vskip 0.2in
 \renewcommand{\theequation}{\thesection.\arabic{equation}}
\section*{2. Statement of results}
\setcounter{section}{2}\setcounter{equation}{0}

Let $\Omega_1$ and $\Omega_2$ be convex open neighbourhoods of the origin in $\mathbb{R}^{n_{11}}$ and $\mathbb{R}^{n_{12}}$, respectively, $\Omega=\Omega_1\times \Omega_2$, the parameter set $\Pi$ be a convex bounded open set of positive Lebesgue measure in $\mathbb{R}^{n_{3}}$. Let $|x|$ denote the maximum norm and $|x|_p$ the $p$-norm. In the following, $l$ and $\alpha$ represent the differentiability orders of functions in the space variables $(I,\varphi)$ and the parameter variables $\xi$, respectively.

\begin{definition}\label{def1}
Let $\alpha$ be a positive integer and $l>0$, $C^{l,\alpha}(\Omega\times \mathbb{T}^{n_2}, \Pi)$ be the class of all functions $f$ on $\Omega\times \mathbb{T}^{n_2}\times \Pi$ whose partial derivatives $\partial_{\xi}^{\beta}f$ with respect to the parameter variable $\xi\in\Pi$ (which means the Whitney derivative if $\Pi$ is a closed set) for all $\beta, 0\leq |\beta|_1\leq \alpha$ are of class $C^l$ in the space variable $x=(I,\varphi)\in  \Omega\times \mathbb{T}^{n_2}$, that is, there is some positive constant $M$ such that the partial derivatives $D^k\left(\partial_{\xi}^{\beta}f\right)$ of $\partial_{\xi}^{\beta}f $ with respect to the space variable $x=(I,\varphi)\in  \Omega\times \mathbb{T}^{n_2}$ satisfy
\begin{equation}\label{1.2}
\left|D^k\left(\partial_{\xi}^{\beta}f(x,\xi)\right)\right|\leq M
\end{equation}
and
\begin{equation}\label{1.3}
\left|D^k\left(\partial_{\xi}^{\beta}f(x,\xi)\right)-D^k\left(\partial_{\xi}^{\beta}f(y,\xi)\right)\right|\leq M|x-y|^{l-[l]}, \quad |k|_1=[l]
\end{equation}
for all $x,y\in \Omega\times \mathbb{T}^{n_2}$ and all $\beta,k$ with $0\leq |\beta|_1\leq \alpha, 0\leq |k|_1\leq [l]$, where $[l]$ is the integer part of $l:\, l-[l]\in [0,1)$, for nonnegative integer vectors $k, \beta$, $D^k=D_1^{k_1}\circ D_2^{k_2}\circ \cdots \circ D_{n_1+n_2}^{k_{n_1+n_2}}$, $D_j^{k_j}=\frac{\partial^{k_j}}{\partial x_j^{k_j}}$, $\partial^\beta_\xi=\frac{\partial^{|\beta|_1}}{\partial \xi_1^{\beta_1}\cdots \partial \xi_{n_3}^{\beta_{n_3}}}$.

In addition, define a norm
$$||f||_{l,\alpha;\Omega\times \mathbb{T}^{n_2},\Pi}=\inf M,$$
which is the smallest M for which the inequalities \eqref{1.2} and \eqref{1.3} hold. $C^{l,\alpha}(\Omega\times \mathbb{T}^{n_2}, \Pi)$ is a Banach space with respect to the norm $||\cdot||_{l,\alpha;\Omega\times \mathbb{T}^{n_2},\Pi}$, which is a generalization of the H\"{o}lder space to a parameter-depending case.  The norms $||\cdot||_{l,\alpha;\mathbb{T}^{n_2},\Pi}$ and $||\cdot||_{\alpha;\Pi}$ are defined in a similar way, which means that the associated function only depends on $\varphi\in \mathbb{T}^{n_2}, \xi\in \Pi$ and $\xi\in \Pi$, respectively.

When $l$ is integer, we also introduce a generalization of the Zygmund space $\hat{C}^{l,\alpha}(\Omega\times \mathbb{T}^{n_2}, \Pi)$ of all functions satisfying
\begin{equation}\label{1.4}
\left|D^k\left(\partial_{\xi}^{\beta}f(x,\xi)\right)\right|\leq M, \qquad 0\leq |k|_1\leq l-1
\end{equation}
and
\begin{equation}\label{1.5}
\left|D^k\left(\partial_{\xi}^{\beta}f(x,\xi)\right)+D^k\left(\partial_{\xi}^{\beta}f(y,\xi)\right)-2D^k\left(\partial_{\xi}^{\beta}f(\frac{1}{2}(x+y),\xi)\right)\right|\leq M|x-y|,  \quad |k|_1=l-1,
\end{equation}
instead of \eqref{1.2} and \eqref{1.3}, respectively, and the norm $||f||_{\hat{C}^l,\alpha;\Omega\times \mathbb{T}^{n_2},\Pi}$
 is the smallest $M$ for which the inequalities \eqref{1.4} and \eqref{1.5} hold. For non-integer $l>0$, $\hat{C}^{l,\alpha}(\Omega\times \mathbb{T}^{n_2}, \Pi)=C^{l,\alpha}(\Omega\times \mathbb{T}^{n_2}, \Pi)$.
\end{definition}

 Sometimes we shall drop parameters from functions whenever there is no confusion.

\vskip 0.2in

 a) Assume

{\bf (H1)} \,  these non-negative constants $q_1, \cdots, q_7$ satisfy
$$q_1>q_3\geq q_5,\quad q_7\geq q_2+q_5,\quad 0<q_2\leq \min\{q_4,q_6\};$$

\vskip 0.2in
{\bf (H2)} \, $\omega_i,A_i\in C^{\alpha}(\Pi)$ with some positive integer $\alpha$ and $A_i$ is a diagonalizable matrix, $A_i(\xi,\varepsilon) =B_i(\xi,\varepsilon)\Lambda_i(\xi,\varepsilon)B_i(\xi,\varepsilon)^{-1}$ for some diagonal matrix $\Lambda_i, i=1,2$.
 Denote $\omega(\xi,\varepsilon)= {\rm col}(\varepsilon^{q_5}\omega_1, \omega_2), \,\Lambda_1(\xi,\varepsilon)={\rm diag} (\lambda_1, \cdots, \lambda_{n_{11}})$, $\Lambda_2(\xi,\varepsilon)= {\rm diag} (\lambda_{n_{11}+1}, \cdots, \lambda_{n_{1}})$. Furthermore, assume that there are positive constants $c_0, c_1$ and $\varepsilon*$ such that for all $\varepsilon \in (0,\varepsilon*]$
$$\inf_{\xi\in \Pi}|\lambda_j(\xi,\varepsilon)|\geq c_0,\quad \inf_{\xi\in \Pi}|\lambda_j(\xi,\varepsilon)-\lambda_i(\xi,\varepsilon)|\geq c_0,\quad i\neq j,\, 1\leq i,j\leq n_{11},\,{\rm or} \, n_{11}+1\leq i,j\leq n_1,$$
$$||B_i||_{\alpha;\Pi},\, ||B_i^{-1}||_{\alpha;\Pi},\, ||\Lambda_i||_{\alpha;\Pi}\leq c_1,\quad ||\omega_i||_{0;\Pi}\coloneqq \sup_{\xi\in \Pi}|\omega_i| \leq c_1,\quad i=1,2,$$
$$\parallel \partial_{\xi}^{\beta}\omega \parallel_{\Pi}\coloneqq \sup_{\xi\in\Pi}| \partial_{\xi}^{\beta}\omega(\xi,\varepsilon) |\leq c_1 \varepsilon^{q_5}, \quad 1\leq |\beta|_1\leq \alpha;$$

\vskip 0.2in
{\bf (H3)} \, $g_j\in C^{l,\alpha}(\Omega\times \mathbb{T}^{n_2}, \Pi)(j=1,\cdots,4)$ with $l>2(\alpha+1)(\iota+2)+\alpha \iota,\, \iota>\alpha n_2-1$.

\vskip 0.1in
{\bf Remark 2.1} The requirement that $A_i$ does not have multiple eigenvalues is not necessary, only for the sake of simplification. The difficulty caused by multiple eigenvalues may be overcome by the technique of R\"{u}ssmann \cite{Russ01}.

\begin{theorem}\label{theorem1}
Suppose that the system \eqref{1.1} satisfies Assumptions (H1)-(H3). Then for any given $0<\gamma\ll1$, there is a sufficiently small $0<\varepsilon_0^*=o(\gamma^{\frac{\alpha+1}{q_2}})$ such that for $0<\varepsilon\leq \varepsilon_0^*$, there exists a Cantor set $\Pi_{\gamma}\subset \Pi$ and for each $\xi \in \Pi_{\gamma}$, the system \eqref{1.1} admits a quasi-periodic invariant torus of the form $I_1=\Phi_1(\varphi;\xi),\, I_2=\Phi_2(\varphi;\xi),\, \varphi={\rm col}(\varphi_1,\varphi_2) \in \mathbb{T}^{n_{21}}\times \mathbb{T}^{n_{22}}$ with frequencies $\omega^*(\xi)=(\varepsilon^{q_5}\omega_1^*(\xi), \omega_2^*(\xi))$, which
is of class $C^\mu$ in $\xi\in\Pi_{\gamma}$ in the sense of whitney and of class $\hat{C}^{l_1}(l_1=l-(\alpha+\mu+2)(\iota+1)-\alpha-2)$ in $\varphi\in \mathbb{T}^{n_2}$ together with derivatives up to order $\mu-1$ with respect to $\xi$ for $0< \mu\leq \alpha$ (positive integer $\mu$), the frequency map $\omega^*(\xi)$ is of class $C^\alpha $ in $\xi\in\Pi_{\gamma}$ in the sense of Whitney and satisfies
\begin{equation}\label{torusest}
||(\Phi_1,\Phi_2)||_{\hat{C}^{l_1},\mu-1; \mathbb{T}^{n_2},\Pi_\gamma}\leq C\varepsilon^{q_2}\gamma^{-(\mu+1)},
\end{equation}
\begin{equation}\label{freqest}
 ||\omega_1^*-\omega_1||_{C_W^\alpha(\Pi_{\gamma})}\leq C \varepsilon^{q_6},\qquad ||\omega_2^*-\omega_2||_{C_W^\alpha(\Pi_{\gamma})}\leq C \varepsilon^{q_7}.
\end{equation}
Moreover, there exist closed subsets $\Pi_\nu$ of $\Pi$, frequency vectors $\omega^\nu(\xi)={\rm col}(\varepsilon^{q_5}\omega_1^\nu(\xi), \omega_2^\nu(\xi))$ and diagonal matrices $\Lambda^\nu(\xi)={\rm diag}(\varepsilon^{q_1}\Lambda_1^\nu(\xi), \varepsilon^{q_3}\Lambda_2^\nu(\xi))$ for $\nu=1,2,\cdots$, satisfying
\begin{equation}\label{frees}
||\omega_1^\nu-\omega_1||_{\alpha;\Pi_{\nu}}\leq c\varepsilon^{q_6},\quad ||\omega_2^\nu-\omega_2||_{\alpha;\Pi_{\nu}}\leq c\varepsilon^{q_7},\quad
||\Lambda_1^\nu-\Lambda_1||_{\alpha;\Pi_{\nu}}\leq c\varepsilon^{q_2},\quad ||\Lambda_2^\nu-\Lambda_2||_{\alpha;\Pi_{\nu}}\leq c\varepsilon^{q_4}
\end{equation}
and
$$\Pi_\nu=\Pi_{\nu-1}\setminus\bigcup_{k,m}\mathfrak{R}_{km}^\nu(\gamma)$$
such that $\Pi_{\gamma}=\bigcap_{\nu=0}^{\infty}\Pi_\nu$, where
$$\mathfrak{R}_{km}^\nu(\gamma)=\left\{\xi\in \Pi_{\nu-1}:\, |\sqrt{-1}\langle k,\omega^{\nu-1}\rangle+\langle m,\Lambda^{\nu-1}\rangle|<\gamma \varepsilon^{q_5}|k|_2^{-\iota}\right\}$$
for $m\in\mathfrak{m}, k\in \mathbb{Z}^{n_2}$,
 $K_{\nu-1}<|k|_2\leq K_\nu$, $\mathfrak{m}=\{m\in \mathbb{Z}^{n_1}:\,|m|_1\leq 2, \sum_{j=1}^{n_1}m_j=0\, {\rm or}\, -1\}$, $\omega^0=\omega, \Lambda^0=\Lambda\coloneqq {\rm diag}(\varepsilon^{q_1}\Lambda_1, \varepsilon^{q_3}\Lambda_2), K_0=0$ and $K_\nu=[K_\nu^\prime ]+1, K_\nu^\prime =r^{-1}3^{\nu+1}[(l+(n_2+1)(\nu+1)-\alpha)\ln 3+(n_2+1)|\ln r|+\ln \widetilde{C}]$, $\Pi_0$ is a closed subset of $\Pi$ whose distance to the boundary of $\Pi$ is at least equal to $\gamma$, where $C$ and $c$ are constants independent of $\nu$, $\varepsilon$ and $\gamma$, $r$ is the radius of the neighbourhood $\Omega$, $||\cdot ||_{C_W^\alpha(\Pi_{\gamma})}$  is the whitney norm (see Appendix A.1), $ \widetilde{C}=24(n_2!)n_2^{n_2}e^{-n_2}$, $[K_\nu^\prime ]$ represents the integer part of $K_\nu^\prime$.
\end{theorem}

Here, we drop $\varepsilon$ from functions, the continuous differentiability of functions $\omega_i^\nu$ and $\Lambda_i^\nu (i=1,2)$ on the  closed set $\Pi_\nu$ means that they are continuously differentiable in some neighbourhood of $\Pi\nu$. Here and in the sequel, we also regard the $\Lambda$ as a column vector of its diagonal elements when $\Lambda$ is a diagonal matrix.

\vskip 0.1in
b) The Cantor set $\Pi_{\gamma}$ is not empty and indeed the measure ${\rm meas}(\Pi\setminus \Pi_{\gamma})\rightarrow 0$ as $\gamma \rightarrow 0$ as long as we impose proper non-degeneracy conditions on frequencies. Since in applications the non-degeneracy conditions on frequencies are different, Theorem 1 does not involve the measure estimate of $\Pi_{\gamma}$ so that it can be used more widely. In the following theorem, we give some conditions to ensure that the Cantor set $\Pi_{\gamma}$ is not empty.

By the assumption (H2), we can write $\omega_2$ as
$$\omega_2(\xi,\varepsilon)=\omega_{20}+\varepsilon^{q_5} \omega_{21}(\xi)+o(\varepsilon^{q_5}),$$
where $\omega_{20}$ is independent of $\xi$,  $o(\varepsilon^{q_5})$ represents infinitely small quantity of $\varepsilon^{q_5}$ up to $\alpha$-th derivatives. Denote $\widetilde{\Lambda}_2(\xi)=\Lambda_2(\xi, \varepsilon)|_{\varepsilon=0}$ and $\widetilde{\omega}(\xi)={\rm col}(\omega_1(\xi,0), \omega_{21}(\xi))$ in the case $n_{22}\neq 0$, $\widetilde{\omega}(\xi)=\omega_1(\xi,0)$ in the case $n_{22}=0$.

\begin{theorem}\label{theorem2}
Suppose that the system \eqref{1.1} satisfies the assumptions in Theorem 1, moreover assume that $n_3+\cdots+n_3^{\alpha}\geq n_2$ and\\

{\rm (i)}\,  for all $\xi\in \Pi$
\begin{equation}\label{I11}
{\rm rank}\left(\widetilde{\omega},\frac{\partial^{|\beta|_1}\widetilde{\omega}}{\partial\xi^{\beta}}:\, 1\leq |\beta|_1\leq \alpha\right)=n_2 \quad {\rm in\, Case }\,\,n_{22}=0,
\end{equation}
\begin{equation}\label{I12}
{\rm rank}\left(\frac{\partial^{|\beta|_1}\widetilde{\omega}}{\partial\xi^{\beta}}:\, 1\leq |\beta|_1\leq \alpha\right)=n_2 \quad {\rm in\, Case }\,\,n_{22}\neq 0,
\end{equation}

{\rm (ii)}\,  for all integer vectors $0\neq k\in \mathbb{Z}^{n_2}, m=(m_1,\cdots,m_{n_{12}})\in \mathbb{Z}^{n_{12}}$ with $1\leq |m|_1\leq 2$ and $m_1+\cdots+m_{n_{12}}=0$ or $-1$
\begin{equation}\label{I2}
{\rm meas}\left\{\xi\in\Pi:\, \sqrt{-1}\langle k,\widetilde{\omega}_0+\varepsilon^{q_5}\widetilde{\omega}(\xi)\rangle+ \varepsilon^{q_3} \langle m,\widetilde{\Lambda}_2(\xi)\rangle =0  \right\}=0,
\end{equation}
where $\widetilde{\omega}_0={\rm col}(0, \omega_{20})$.

Then the Cantor set $\Pi_{\gamma}$ defined in Theorem 1 is of positive Lebesgue measure and ${\rm meas}(\Pi\backslash\Pi_{\gamma})\rightarrow 0$ as $\gamma \rightarrow 0$.
\end{theorem}

The conditions \eqref{I11}, \eqref{I12}, Lemma~\ref{lemmaA.3} and Remark A.2 imply that there is a constant $c_2>0$ such that
\begin{equation}\label{Dn11}
\max_{0\leq \mu \leq \alpha}\parallel D^{\mu}\langle b,\widetilde{\omega}(\xi)\rangle\parallel \geq c_2 \qquad {\rm in \, Case}\,\, n_{22}=0
\end{equation}
and
\begin{equation}\label{Dn22}
\max_{1\leq \mu \leq \alpha}\parallel D^{\mu}\langle b,\widetilde{\omega}(\xi)\rangle\parallel \geq c_2 \qquad {\rm in \, Case}\,\, n_{22}\neq 0
\end{equation}
for all $\xi\in \Pi, b\in \mathcal{S}_{n_2,1}\coloneqq \{b\in \mathbb{R}^{n_2}:\, |b|_2=1\}$.\\
Let
$$K^*=\frac{32c_1}{c_2}n_3^{\frac{\alpha}{2}},\qquad f_{km}(\xi)=\langle k,\widetilde{\omega}_0+\varepsilon^{q_5}\widetilde{\omega}(\xi)\rangle+ \varepsilon^{q_3} \langle m,{\rm Im}\widetilde{\Lambda}_2(\xi)\rangle,$$
where ${\rm Im}\widetilde{\Lambda}_2$ is the imaginary part of $\widetilde{\Lambda}_2$.

\begin{theorem}\label{theorem3}
If in Theorem 2, the condition (ii) is replaced by

{\rm (ii)'}\,  there is a constant $c_3>0$ such that
$$\max_{1\leq \mu \leq \alpha}\parallel D^{\mu}f_{km}(\xi)\parallel \geq c_3\varepsilon^{q_5} \qquad {\rm for \, all}\,\, \xi\in\Pi,0<|k|_2<K^*,$$
then
$${\rm meas}\Pi_{\gamma}={\rm meas}\Pi-O(\gamma^{\frac{1}{\alpha}})$$
for sufficiently small $\gamma$.
\end{theorem}

{\bf Remark 2.2} If $q_3>q_5\geq 0$, or there a constant $c_2^{\prime}>0$ such that
 $$\inf_{\xi\in\Pi}|\langle m,{\rm Re}\widetilde{\Lambda}_2(\xi)\rangle|\geq c_2^{\prime} \qquad {\rm for}\,1\leq |m|_1\leq 2,\,m_1+\cdots+m_{n_{12}}=0\,{\rm or}\, -1,$$
then the conditions {\it (ii)} and {\it (ii)'} in Theorems 2 and 3, respectively, may be removed, see the proof of Theorems 2 and 3,  and Remark 5.1 in Section 5.

More results on measure estimates of $\Pi_{\gamma}$ will be given in the forthcoming second part concerning on the persistence of quasi-periodic invariant tori in bifurcation theory.

\vskip 0.1in
c) we consider a specific form of \eqref{1.1} for the case $n_{11}=n_{22}=0, n_2=n_3, q_3=q_5=0$ and $q_4=q_6=1$, which means that the first and fourth equations in \eqref{1.1} are absent and the number of parameter variables equals the dimension of tori, the equation \eqref{1.1} reads
\begin{equation}\label{eq1}
\left\{\begin{array}{rl}
\dot I & = A(\xi)I+\varepsilon g_1(I,\varphi; \xi,\varepsilon)\\
\dot{\varphi} & =\omega(\xi)+\varepsilon g_2(I,\varphi; \xi,\varepsilon).
\end{array}
\right.
 \end{equation}
 Denote $\Lambda={\rm diag} (\lambda_1, \cdots, \lambda_{n_{1}})$ and $\omega={\rm col}(\omega_1,\cdots, \omega_{n_2})$, where $\lambda_1, \cdots, \lambda_{n_{1}} $ are the eigenvalues of $A$, $A(\xi) =B(\xi)\Lambda(\xi)B(\xi)^{-1}$. Assume
 \vskip 0.2in
{\bf (H2)$^{\prime}$} \, $\omega,A\in C^1(\Pi)$, the map $\xi\rightarrow \omega(\xi)$ is a diffeomorphism between $\Pi$ and its image, and there exist positive constants $c_0, c_1$ and $c_4$ such that $||B||_{1;\Pi}, ||B^{-1}||_{1;\Pi}, ||\Lambda||_{1;\Pi}, ||\omega||_{1;\Pi}\leq c_1$,
\begin{equation}\label{normalfre}
|\langle m,\Lambda(\xi)\rangle| \geq c_0, \quad \left|\left| \left(\frac{\partial\omega}{\partial \xi}\right)^{-1}\right|\right|\leq c_4 \qquad {\rm on}\quad \Pi
\end{equation}
and
\begin{equation}\label{mixfre}
{\rm meas}\left\{\xi\in\Pi:\, \sqrt{-1}\langle k,\omega(\xi)\rangle+  \langle m,\Lambda(\xi)\rangle =0  \right\}=0
\end{equation}
for all integer vectors $0\neq k\in \mathbb{Z}^{n_2}, m\in \mathbb{Z}^{n_1}$ with $1\leq |m|_1\leq 2$ and $m_1+\cdots+m_{n_1}=0$ or $-1$;
\vskip 0.2in
{\bf (H3)$^{\prime}$} \, $g_j\in C^{l,1}(\Omega\times \mathbb{T}^{n_2}, \Pi)(j=1,2)$ with $l>5\iota+8,\, \iota>n_2-1$.

\vskip 0.2in
{\bf Remark 2.3} \,(i) When the real part ${\rm Re} \Lambda$ of $\Lambda$ satisfies $\langle m,{\rm Re}\Lambda(\xi)\rangle \neq 0$ on $\Pi$, the condition \eqref{mixfre} holds spontaneously. In particular, the condition \eqref{mixfre} is satisfied if $\Lambda$ is independent of $\xi$.

(ii)  The Assumption (H2)$^{\prime}$ implies that the condition \eqref{mixfre} is satisfied if
\begin{equation}\label{Dfre}
\left(\left(\frac{\partial\omega}{\partial \xi}\right)^{-1}\right)^T\frac{\partial}{\partial \xi}\langle m,\Lambda(\xi)\rangle\neq \sqrt{-1} k \quad {\rm for}\, 0\neq |k|_1\leq 2n_2 c_1 c_4.
\end{equation}

Theorems 1-3 imply
\begin{corollary}\label{corollary1}
Suppose that the system \eqref{eq1} satisfies Assumptions (H2)$^{\prime}$ and (H3)$^{\prime}$. Then for any given $0<\gamma \ll 1$, there is a sufficiently small $\varepsilon^*>0$ such that for $0<\varepsilon\leq \varepsilon^*$, there exists a Cantor set $\Pi_{\gamma}\subset \Pi$ with positive Lebesgue measure (the measure satisfies the estimate
${\rm meas}\Pi_{\gamma}={\rm meas}\Pi-c\gamma$ if \eqref{Dfre} replaces \eqref{mixfre})
and for each $\xi\in \Pi_{\gamma}$, the system \eqref{eq1} possesses a quasi-periodic invariant torus $I=\Phi(\varphi;\xi),\, \varphi \in \mathbb{T}^{n_2}$ consisting of quasi-periodic motions, which is of $\hat{C}^{l_1}(l_1=l-4(\iota+1)-3)$ in $\varphi \in \mathbb{T}^{n_2}$ and Lipschitz in  $\xi\in \Pi_{\gamma}$, where $c$ is a constant independent of $\gamma$ and $\varepsilon$.
\end{corollary}

Usually the normal form (integrable part) of \eqref{1.1} related bifurcation problems of actual models is only finitely differentiable, not analytic in the parameter $\xi$, and the frequency map is possibly degenerate so that we need the higher-order derivatives of the frequency map to estimate the Lebesgue measure of $\Pi_{\gamma}$ and obtain $\Pi_{\gamma}$ is the most part of $\Pi$. Hence, we want to establish an approximation lemma and the corresponding inverse approximation lemma in which a finitely differentiable function is approximated by a sequence of functions being analytic in space variables, but finitely differentiable in parameter variables. These comprise Section 3. The proofs of Theorems 1-3 are given in Sections 4 and 5, respectively.

\vskip 0.2in
\renewcommand{\theequation}{\thesection.\arabic{equation}}
\section*{3  Approximation Lemmas}
\setcounter{section}{3}\setcounter{equation}{0}

Zehnder \cite{Zehn75} established the approximation and inverse approximation Lemmas on a finitely differentiable real function approximated by a sequence of real analytic functions, which was generalized to the anisotropic case by P\"{o}schel \cite{Pos82}, and was sharpened to covering the finitely differentiable and Gevrey regular cases by Wagener \cite{Wag10}, respectively. Here, we give generalized versions of Zehnder's approximation and inverse approximation lemmas finite- smoothly depending on parameters, and obtain estimates of higher-order regularity.
\vskip 0.2in

a) We first introduce some notations. Let $m,n$ and $\alpha$ be positive integers, $\mathcal{U}\subset \mathbb{C}^m$ and $\Pi\subset \mathbb{R}^n$ be open sets, $\mathfrak{A}^\alpha (\mathcal{U},\Pi)$ be the class of all functions of $(z,\xi)$ on $\mathcal{U}\times\Pi$ which are analytic in $z\in \mathcal{U}$ and $\alpha$-times continuously differentiable in $\xi\in \Pi$. For $g\in \mathfrak{A}^\alpha (\mathcal{U},\Pi)$, define
$$|g|_{\mathcal{U},\alpha;\Pi}=\sup_{|\beta|_1\leq \alpha}\sup_{(z,\xi)\in \mathcal{U}\times\Pi}\left|\partial_\xi^\beta g(z,\xi)\right|.$$
In particular, for $\mathcal{U}=\{z\in \mathbb{C}^m:\, |{\rm Im} z|\coloneqq sup_{1\leq j\leq m}|{\rm Im} z_j|<r\}$, we denote $|g|_{\mathcal{U},\alpha;\Pi}$ by $|g|_{r,\alpha;\Pi}$.

Take an even function $u_0\in C_0^\infty(\mathbb{R})$, vanishing outside the interval $[-1,1]$ and identically equal to 1 in a neighbourhood of 0 (see \cite{Wag10} for the construction of such a function). For $x\in \mathbb{R}^m$, let $u(x)=u_0(|x|_2^2)$ and $\tilde{u}$ be the inverse Fourier transform of $u$
$$\tilde{u}(z)=(2\pi)^{-m}\int_{\mathbb{R}^m} u(x)e^{\sqrt{-1}\langle z,x\rangle}dx.$$
Let $f$ be a real-valued function of class $C^{l,\alpha}(\mathbb{R}^m, \Pi)$  (see Definition\ref{def1}), $f_r\,(0<r\leq 1)$ be defined by the convolution
\begin{equation}\label{frexpr}
f_r(x,\xi)\coloneqq (S_r(f(\cdot,\xi))(x)=r^{-m}\int_{\mathbb{R}^m} \tilde u(r^{-1}(x-y))f(y,\xi)dy
\end{equation}
for $x\in \mathbb{C}^m$. We list some properties of the analytic smoothing operator $S_r$ in Section A.3 of the appendix, which will be used in the proof of the next lemma.

\begin{lemma}\label{lemma2.1}
 Let $f(x,\xi)$ be a real-valued function of class $C^{l,\alpha}(\mathbb{R}^m,\Pi)$ for some real number $l> 0$ and $\alpha \in \mathbb{N}$, where $\Pi\subset \mathbb{R}^n$ is an open set. Then
for every $r\in (0,1]$, the function $f_r(x,\xi)$ is $\alpha$-times continuously differentiable in $\xi\in \Pi$, entire real analytic in $x\in \mathbb{C}^m$ together with derivatives up to order $\alpha$ with respect to $\xi$, and satisfies

{\rm (i)}\, $||f_r-f||_{p,\alpha;\mathbb{R}^m,\Pi}\leq C_1(l,p)r^{l-p}||f||_{l,\alpha;\mathbb{R}^m,\Pi}$\,\, for all $0\leq p\leq l$,

{\rm (ii)}\, $|f_r-f_{r^\prime}|_{r^\prime,\alpha;\Pi}\leq C_2(l,p)r^{p}||f||_{p,\alpha;\mathbb{R}^m,\Pi}$\,\, for all $0\leq p\leq l$ and $0<r^\prime\leq r$,

{\rm (iii)}\,  $|f_r|_{r,\alpha;\Pi}\leq C_3(l)||f||_{0,\alpha;\mathbb{R}^m,\Pi}\leq C_3(l)||f||_{l,\alpha;\mathbb{R}^m,\Pi}$,\\
where $C_j\,(j=1,2,3)$ are constants depending on $l,p$ and the dimension $m$. Moreover, $f_r$ is $\omega$-periodic in some variable if in which $f$ is $\omega$-periodic.
\end{lemma}

{\bf Proof}  From \eqref{frexpr} it is clear that $f_r(x,\xi)$ is  analytic in $x\in \mathbb{C}^m$, and $\alpha$-times continuously differentiable in $\xi\in \Pi$, taking real values on real variables $x$, and if $f$ is periodic in some variable, then so is $f_r$. As differentiation may commute with integration in \eqref{frexpr} for functions with bounded derivatives, we obtain $\partial_\xi^\beta f_r=S_r(\partial_\xi^\beta f)$ for $|\beta|_1\leq \alpha$. Of course, we also have $S_r(D^kf)=D^k(S_r(f))$ for $|k|_1\leq l, k\in \mathbb{Z}^m$. Hence we only need to prove the estimates (i)-(iii) in the case without parameter-dependence. In the following, we will use $C$ to denote some constant depending $l, p$ and $m$.

(i) The case where $p$ is a integer, is proved by Chierchia \cite{Chier03}, see Lemma~\ref{lemmaA4}(f) in Appendix. Hence we only give the proof for the case $p=q+\mu\leq l, \mu\in (0,1), q\in \mathbb{Z}_+$. Denote $g(x)=D^\beta f, |\beta|_1=q$. Then by (a) and (b) in Lemma~\ref{lemmaA4}, we have for $x,y\in \mathbb{R}^m$,
\begin{eqnarray*}
& &  \sup_{x\neq y}|x-y|^{-\mu}|(g-S_rg)(x)-(g-S_rg)(y)|\\
\hskip 0.2in & & =\sup_{x\neq y}|x-y|^{-\mu}\left|\int_{\mathbb{R}^m}\tilde{u}(z)[g(x)-g(x-rz)-g(y)+g(y-rz)]dz\right|\equiv (*).
\end{eqnarray*}

Case I: $q=[l]$, the integer part of $l$. For $|x-y|\geq r$, by $g\in C^{l-q}$ and Lemma~\ref{lemmaA4} (d), we have
\begin{eqnarray*}
 (*) & \leq & \sup_{x\neq y}|x-y|^{-\mu}\int_{\mathbb{R}^m}|\tilde{u}(z)|(|g(x)-g(x-rz)|+|g(y)-g(y-rz)|)dz\\
 & \leq & 2r^{l-p}||f||_{l;\mathbb{R}^m}\int_{\mathbb{R}^m}|\tilde{u}(z)||z|^{l-q}dz\leq C r^{l-p}||f||_{l;\mathbb{R}^m}.
\end{eqnarray*}

For $|x-y|< r$, we also have
\begin{eqnarray*}
 (*) & \leq & \sup_{x\neq y}|x-y|^{-\mu}\left(||f||_{l;\mathbb{R}^m}|x-y|^{l-q} + \int_{\mathbb{R}^m}|\tilde{u}(z)||g(x-rz)-g(y-rz)|dz\right)\\
 & \leq & \left(1+\int_{\mathbb{R}^m}|\tilde{u}(z)|dz\right)|x-y|^{l-p}||f||_{l;\mathbb{R}^m}\leq C r^{l-p}||f||_{l;\mathbb{R}^m}.
\end{eqnarray*}
Hence, $||g-S_rg||_{\mu;\mathbb{R}^m}\leq C r^{l-p}||f||_{l;\mathbb{R}^m}$, which, combining with Lemma~\ref{lemmaA4} (f) for the case of integers, implies (i) for the case $q=[l]$.

Case II: $q<[l]$. For $|x-y|\geq r$, using the Taylor's formula of $h(rz)=g(x-rz)-g(y-rz)$ at $z=0$ and Lemma~\ref{lemmaA4} (c), we obtain
\begin{eqnarray}
(*) &= &  \sup_{x\neq y}\frac{|x-y|^{-\mu}}{([l]-q)!}
\left|\int_{\mathbb{R}^m}\tilde{u}(z)(-rz\cdot \nabla )^{[l]-q}(g(x-\theta rz)-g(y-\theta rz))dz\right| \nonumber\\
& \leq &  \sup_{x\neq y}|x-y|^{-\mu}r^{[l]-q}\sum_{|k|_1=[l]-q}\frac{1}{k!}
\left|\int_{\mathbb{R}^m}\tilde{u}(z)z^k (D^kg(x-\theta rz)-D^kg(y-\theta rz))dz\right|,\label{cpest1}
\end{eqnarray}
where $rz\cdot \nabla=\sum_{j=1}^{m}rz_jD_j$ and $\theta\in (0,1)$. Thus, equivalently, we need to estimate the following expression
$$(**)\equiv  \sup_{x\neq y}|x-y|^{-\mu}
\left|\int_{\mathbb{R}^m}\tilde{u}(z)z^k (D^{k+\beta}f(x-\theta rz)-D^{k+\beta}f(y-\theta rz))dz\right|,\, |k+\beta|_1=[l].$$
By Lemma~\ref{lemmaA4} (c), we get
\begin{eqnarray}
(**) &= &  \sup_{x\neq y}|x-y|^{-\mu}\left|\int_{\mathbb{R}^m}\tilde{u}(z)z^k (D^{k+\beta}f(x-\theta rz)-D^{k+\beta}f(x))dz\right. \nonumber\\
& &+ \left.\int_{\mathbb{R}^m}\tilde{u}(z)z^k (D^{k+\beta}f(y)-D^{k+\beta}f(y-\theta rz))dz\right|\leq C r^{l-[l]-\mu}||f||_{l;\mathbb{R}^m}.\label{cpest2}
\end{eqnarray}

For $|x-y|< r$, if $[l]-q\geq 2$, then similarly we have
\begin{equation}\label{cpest3}
(*)\leq \sup_{x\neq y}|x-y|^{-\mu}r^{[l]-q-1}\sum_{|k|_1=[l]-q-1}\frac{1}{k!}
\left|\int_{\mathbb{R}^m}\tilde{u}(z)z^k (D^kg(x-\theta rz)-D^kg(y-\theta rz))dz\right|.
\end{equation}
The mean value theorem and Lemma~\ref{lemmaA4} (c) deduce
\begin{eqnarray}
& & \sup_{x\neq y}|x-y|^{-\mu}
\left|\int_{\mathbb{R}^m}\tilde{u}(z)z^k (D^kg(x-\theta rz)-D^kg(y-\theta rz))dz\right| \nonumber\\
& \leq & \sup_{x\neq y}|x-y|^{-\mu}\sum_{|k'|_1=1}
\left|\int_{\mathbb{R}^m}\tilde{u}(z)z^k (x-y)^{k'}D^{k+k'}g(y-\theta rz+\theta_{kk'}(x-y))dz\right| \nonumber\\
& =& \sup_{x\neq y}|x-y|^{-\mu}\sum_{|k'|_1=1}
\left|\int_{\mathbb{R}^m}\tilde{u}(z)z^k (x-y)^{k'}(D^{k+k'}g(y-\theta rz+\theta_{kk'}(x-y))-D^{k+k'}g(x))dz\right| \nonumber\\
&\leq & \sup_{x\neq y}|x-y|^{-\mu}\sum_{|k'|_1=1}
|x-y|||f||_{l;\mathbb{R}^m}\int_{\mathbb{R}^m}|\tilde{u}(z)z^k| |\theta rz+(1-\theta_{kk'})(x-y)|^{l-[l]}dz \nonumber\\
&\leq & C r^{l+1-[l]-\mu}||f||_{l;\mathbb{R}^m},\label{cpest4}
\end{eqnarray}
where $\theta_{kk'}\in (0,1)$. If $[l]-q=1$, then by the mean value theorem,
\begin{eqnarray}\label{cpest5}
(*) & \leq & \sup_{x\neq y}|x-y|^{-\mu}\sum_{|k'|_1=1}
\left|\int_{\mathbb{R}^m}\tilde{u}(z) (x-y)^{k'}(D^{k'}g(y+\theta_{1k'}(x-y))-D^{k'}g(y-rz+\theta_{2k'}(x-y)))dz\right| \nonumber \\
&\leq & \sup_{x\neq y}|x-y|^{-\mu}\sum_{|k'|_1=1}
|x-y|||f||_{l;\mathbb{R}^m}\int_{\mathbb{R}^m}|\tilde{u}(z)| |rz+(\theta_{1k'}-\theta_{2k'})(x-y)|^{l-[l]}dz \nonumber \\
&\leq & C r^{l+1-[l]-\mu}||f||_{l;\mathbb{R}^m}=C r^{l-p}||f||_{l;\mathbb{R}^m},
\end{eqnarray}
where $\theta_{1k'},\theta_{2k'}\in (0,1)$. Hence, \eqref{cpest1}-\eqref{cpest5} and Lemma~\ref{lemmaA4} (f) imply (i) for the case $q<[l]$.

Obviously, Lemma~\ref{lemmaA4} (f) implies (ii), and the definition of $f_r$ and Lemma~\ref{lemmaA4} (b) and (d) imply (iii). \hskip 0.1in $\blacksquare$
\vskip 0.2in

From Lemma~\ref{lemma2.1}, it follows the approximation lemma.
\begin{lemma}\label{apprlem}
 (Approximation Lemma) Let $f(x,\xi)$ be a real-valued function of class $C^{l,\alpha}(\mathbb{R}^m,\Pi)$ for some real number $l> 0$ and $\alpha \in \mathbb{N}$, where $\Pi$ is an open set, and let $\{ r_j\}_{j=0}^{\infty}$ be a monotonically decreasing sequence of positive numbers with $r_0\leq 1$ and tend to zero. Then there exists a sequence of functions $\{f_j(z,\xi)\}_{j=0}^{\infty}$, being of class $C^{\alpha}$ in $\xi \in \Pi$, and entire, real analytic in $z\in \mathbb{C}^m$ together with derivatives up to order $\alpha$ with respect to $\xi$, starting with $f_0\equiv 0$, such that
$$\lim_{j\rightarrow \infty}||f_j-f||_{p,\alpha;\mathbb{R}^m,\Pi}=0 \qquad {\rm for \,\, all} \,\, 0\leq p<l$$
and
$$|f_j-f_{j-1}|_{r_j,\alpha;\Pi}\leq C_0 r_{j-1}^l||f||_{l,\alpha;\mathbb{R}^m,\Pi} \qquad {\rm for } \,\, j\geq 1,$$
where the constant $C_0$ depends on $l$ and the dimension $m$. Moreover, the $f_j$ is $\omega$-periodic in each variable in which $f$ is $\omega$-periodic.
\end{lemma}

\vskip 0.2in

b) Now, we want to apply the approximation lemma to the proof of Theorem 1 and obtain sequences of real analytic functions approximating $g_i(i=1,\cdots,4)$ in the equation (1.1).

Without loss of generality, we take
$$\Omega=\{I={\rm col}(I_1,I_2)\in \mathbb{R}^{n_1}:\, |I|<3\tilde{r}\}$$
for some constant $0<\tilde{r}\leq 1$. Let
$$\Omega^*=\{I\in \mathbb{R}^{n_1}:\, |I|\leq 2\tilde{r}\},\qquad r_j=\tilde{r}3^{-j},\, j=0,1,2,\cdots.$$
Define complex neighbourhoods $\mathcal{U}_j$ of $\Omega^*\times \mathbb{T}^{n_2}$ for $j=0,1,2,\cdots$ by
$$\mathcal{U}_j=\{(I,\varphi)\in \mathbb{C}^{n_1}\times  \mathbb{C}^{n_2}:\, {\rm dist}(I,\Omega^*)<3r_j, |{\rm Im} \varphi|<3r_j\}\coloneqq \Omega^*\times \mathbb{T}^{n_2}+(3r_j,3r_j).$$
We first expand the definition domain $\Omega\times \mathbb{T}^{n_2}\times \Pi$ of $g_i(i=1,\cdots,4)$ to $\mathbb{R}^{n_1}\times \mathbb{T}^{n_2}\times \Pi$ in the following manner: we multiply $g_i$ by a $C^\infty$-function on $\mathbb{R}^{n_1}$ which identical 1 on $\Omega^*$ and vanishes outside $\Omega$. The obtained function belongs to $C^{l,\alpha}(\mathbb{R}^{n_1}\times \mathbb{T}^{n_2}, \Pi)$ and is equal to $g_i$ on $\Omega^*\times \mathbb{T}^{n_2}\times \Pi$, its norm is bounded by $c_l||g_i||_{l,\alpha;\Omega\times \mathbb{T}^{n_2},\Pi}$, where $c_l$ is a constant depending $l,n_1$ and the chosen $C^\infty$-function. Then by the approximation lemma (Lemma~\ref{apprlem}) we have the following corollary.

\begin{corollary}\label{corollary2}
If the system (1.1) satisfies Assumption (H3), then there exist sequences $\{g_i^j(I,\varphi,\xi)\}_{j=0}^\infty$ $(i=1,\cdots,4)$ of real analytic functions, being of class $C^{\alpha}$ in $\xi \in \Pi$, and entire, real analytic in $(I,\varphi)\in \mathcal{U}_0$, periodic in the variables $\varphi$ with periodic $2\pi$ together with derivatives up to order $\alpha$ with respect to $\xi$, starting with $g_i^0\equiv 0$, such that
$$\lim_{j\rightarrow \infty}||g_i^j-g_i||_{p,\alpha;\Omega^*\times \mathbb{T}^{n_2},\Pi}=0 \qquad {\rm for \,\, all} \,\, 0\leq p<l,$$
$$|g_i^j-g_i^{j-1}|_{\mathcal{U}_{j-1},\alpha;\Pi}\leq C_0 r_{j-1}^l||g_i||_{l,\alpha;\Omega\times \mathbb{T}^{n_2},\Pi} \qquad {\rm for } \,\, j\geq 1, i=1,\cdots,4, $$
where $C_0$ is a constant depending only on $l,n_1,n_2,\tilde{r}$ and $c_l$.
\end{corollary}
\vskip 0.2in

 c) Let $\Omega\subset \mathbb{R}^m$ be an open convex set, and $\Pi_0\subset \mathbb{R}^n$ be a closed set,
$$\mathcal{W}_j=\Omega+r_j,\quad  \Pi_j=\bigcup_{\xi\in \Pi_0}\{\zeta \in \mathbb{R}^n: \, |\zeta-\xi|<s_j\},\qquad j=0,1,2,\cdots,$$
where $r_j=r_0\theta^j,0<\theta<1$ and $\{s_j\}_{j=0}^{\infty}$ is a  monotonically decreasing sequence of positive numbers with $s_0\leq 1$ and tend to zero.\\

\begin{lemma}\label{invapplem} (Inverse Approximation Lemma)
 Let $\{f_j(x,\xi)\}_{j=0}^{\infty}$ be a sequence of functions such that $f_0\equiv 0$, $f_j(x,\xi)$ is of class $C^{\alpha}$ in $\xi \in \Pi_j$,  real analytic in $x\in \mathcal{W}_j$ together with derivatives up to order $\alpha$ with respect to $\xi$, and
\begin{equation}\label{fjest}
|f_j-f_{j-1}|_{\mathcal{W}_j,\alpha;\Pi_j}\leq M r_j^l
\end{equation}
 for every $ j\geq 1$ and some constant $M$. If there exists a constant $c^\prime_0>0$ such that $r_j^l\leq c^\prime_0s_j^\alpha, j=1,2,\cdots$, then there is a unique function $f(x,\xi)$ being of class $C^{\alpha}$ in $\xi \in \Pi_0$ in the sense of Whitney (see Appendix A.1), and of class $\hat{C}^l$ in $x\in \Omega$ together with derivatives up to order $\alpha-1$ with respect to $\xi$ such that
 $$||f||_{\hat{C}^l,\alpha-1;\Omega,\Pi}\leq C^\prime_0 M \qquad {\rm and} \quad
\lim_{j\rightarrow \infty}||f-f_j||_{p,\alpha-1;\Omega,\Pi}=0 \qquad {\rm for \,\, all} \quad 0\leq p<l.$$

  Moreover, let $l=q+\mu,q\in \mathbb{Z}_+, \mu>0$ and if $r_j^\mu\leq c^\prime_1s_j^\delta$ for some constant $c^\prime_1$ and $0<\delta\leq 1$, then we may require the $(\alpha-1)$-order derivatives $\partial_\xi^\beta f(x,\xi)$ with $|\beta|_1= \alpha -1$ to be uniformly $\delta$-H\"{o}lder continuous in $\xi\in \Pi_0$ in the space $C^q(\Omega)$, that is,
 \begin{equation}\label{lipcon}
 ||\partial_\xi^\beta f(\cdot,\xi)-\partial_\xi^\beta f(\cdot,\zeta)||_{C^q(\Omega)}\leq C^\prime_1M|\xi-\zeta|^{\delta} \quad {\rm for}\,\xi,\zeta\in \Pi_0,\,|\beta|_1= \alpha -1,
 \end{equation}
   where the constant $C^\prime_0$ and $C^\prime_1$ depend on $l,m,n,\theta,c^\prime_0$ and $c^\prime_1$, $\hat{C}^l(\Omega)$ is the Zygmund space.
\end{lemma}

{\bf Proof} By a similar proof to that of Lemma 2.2 (ii) in \cite{Zehn75} (also see the proof of Lemma 4.3 in \cite{LL10}, Theorem A.3 in \cite{Wag10}), we can obtain that there exist functions $f^{(\beta)}\in \hat{C}^l(\Omega), |\beta|_1\leq \alpha$ such that
$$\sup_{\xi\in \Pi}||f^{(\beta)}(\cdot,\xi)||_{\hat{C}^l(\Omega)}\leq C^\prime_0 M \quad {\rm and} \quad
\lim_{j\rightarrow \infty}||\partial_\xi^\beta f_j(\cdot,\xi)-f^{(\beta)}(\cdot,\xi)||_{C^p(\Omega)}=0$$
uniformly on $\Pi_0$ for all $0\leq p<l$ and $|\beta|_1\leq \alpha$. Set $f(x,\xi)= f^{(\beta)}(x,\xi)$ with $\beta=0$. To prove the rest of the lemma we only need to verify \eqref{lipcon} and the compatibility conditions in the definition of Whitney derivatives (see Appendix A.1)
\begin{equation}\label{compc}
f^{(\beta)}(x,\xi)=\sum_{|\beta+k|_1\leq \alpha-1}\frac{1}{k!}f^{(\beta+k)}(x,\zeta)(\xi-\zeta)^k+R^\beta(x,\xi,\zeta)
\end{equation}
with
\begin{equation}\label{reterm}
\sup_{x\in\Omega}|R^\beta(x,\xi,\zeta)|\leq C M|\xi-\zeta|^{\alpha-|\beta|_1}
\end{equation}
for all $\xi,\zeta\in \Pi_0, |\beta|_1\leq \alpha-1$ and some finite constant $C$.

Set
$$h_j(x,\xi)=f_j(x,\xi)-f_{j-1}(x,\xi),\quad R_j^\beta(x,\xi,\zeta)=\partial_\xi^\beta h_j(x,\xi)-\sum_{|\beta+k|_1\leq \alpha-1}\frac{1}{k!}
\partial_\xi^{\beta+k} h_j(x,\zeta)(\xi-\zeta)^k$$
for $j\geq 1, |\beta|_1\leq \alpha-1$. Then
\begin{equation}\label{fdex}
f^{(\beta)}(x,\xi)=\sum_{j=1}^{\infty}\partial_\xi^\beta h_j(x,\xi),\quad
R^\beta(x,\xi,\zeta)=\sum_{j=1}^{\infty}R_j^\beta (x,\xi,\zeta),\quad |\beta|_1\leq \alpha-1.
\end{equation}
If $s_{j_0+1}\leq |\xi-\zeta|<s_{j_0}$ for some positive integer $j_0$, then the line segment $L$ connecting $\xi$ to $\zeta$ is contained in $\Pi_j$ with $1\leq j\leq j_0$, and the Taylor expansion implies
$$\sup_{x\in\Omega}|R_j^\beta(x,\xi,\zeta)|\leq C_1(\beta) Mr_j^l|\xi-\zeta|^{\alpha-|\beta|_1},\qquad 1\leq j\leq j_0.$$
And
$$\sup_{x\in\Omega}|R_j^\beta(x,\xi,\zeta)|\leq C_2(\beta) Mr_j^ls_{j_0+1}^{-(\alpha-|\beta|_1)} |\xi-\zeta|^{\alpha-|\beta|_1},\qquad j\geq j_0+1.$$
Hence,
\begin{eqnarray*}
\sup_{x\in\Omega}|R^\beta(x,\xi,\zeta)| &\leq & M  |\xi-\zeta|^{\alpha-|\beta|_1}\left(C_1\sum_{j=1}^{j_0} r_j^l+C_2\sum_{j=j_0+1}^{\infty}\left(\frac{r_j}{r_{j_0+1}}\right)^l\frac{r_{j_0+1}^l}{s_{j_0+1}^{\alpha-|\beta|_1}}\right)\\
& \leq & C M|\xi-\zeta|^{\alpha-|\beta|_1}
\end{eqnarray*}
If $|\xi-\zeta|\geq s_1$, then we also have
\begin{eqnarray*}
\sup_{x\in\Omega}|R^\beta(x,\xi,\zeta)| &\leq & C_2(\beta) M  |\xi-\zeta|^{\alpha-|\beta|_1}\frac{r_1^l}{s_1^{\alpha-|\beta|_1}} \sum_{j=1}^{\infty}\left(\frac{r_j}{r_1}\right)^l\\
& \leq & C M|\xi-\zeta|^{\alpha-|\beta|_1}
\end{eqnarray*}
Thus, we prove the compatibility conditions \eqref{compc} and \eqref{reterm}, and obtain $\partial_\xi^\beta f(x,\xi)=f^{(\beta)}(x,\xi)$ for $|\beta|_1\leq \alpha-1$.

Now, we prove \eqref{lipcon}. Let
$$u_j(x,\xi)=\partial_\xi^\beta h_j(x,\xi) \quad {\rm and}\quad u(x,\xi)=\partial_\xi^\beta f(x,\xi),\qquad |\beta|_1=\alpha-1.  $$
Then the \eqref{fdex} implies
\begin{equation}\label{uexp}
u(x,\xi)=\sum_{j=1}^{\infty}u_j(x,\xi)\qquad {\rm for}\quad (x,\xi)\in \Omega\times \Pi_0.
\end{equation}
By the Cauchy inequality and \eqref{fjest}, we have
\begin{equation}\label{ujder}
|D^k u_j|_{\Omega,1;\Pi_j}\leq C(k) M r_j^{l-|k|_1} \qquad {\rm for} \quad |k|_1\leq q,
\end{equation}
where $C(k)$ is a constant depending only on $k$. By a similar proof to one for the compatibility and replacing \eqref{fjest} with \eqref{ujder}, \eqref{uexp} implies
$$\sup_{x\in \Omega}|D^k u(x,\xi)-D^k u(x,\zeta)|\leq C_1^\prime M|\xi-\zeta|^{\delta} \quad {\rm for}\,\xi,\zeta\in \Pi_0,\,|k|_1\leq q.$$
The proof of the lemma is complete.\hskip 0.4in $\blacksquare$

\vskip 0.2in
\renewcommand{\theequation}{\thesection.\arabic{equation}}
\section*{4  Proof of Theorem 1}
\setcounter{section}{4}\setcounter{equation}{0}

We first introduce some notation so that the system (1.1) is written in a compact form. Denote
$$A^0={\rm diag}(\varepsilon^{q_1}A_1, \varepsilon^{q_3}A_2),\qquad B={\rm diag}(B_1,B_2),$$
$$\Lambda^0={\rm diag}(\varepsilon^{q_1}\Lambda_1, \varepsilon^{q_3}\Lambda_2),\qquad \omega^0={\rm col}(\varepsilon^{q_5}\omega_1,\omega_2),$$
$$P_1={\rm diag}(\varepsilon^{q_1}E_{n_{11}},\varepsilon^{q_3+q_4-q_2}E_{n_{12}}),\quad P_2={\rm diag}(\varepsilon^{q_5+q_6-q_2}E_{n_{21}},\varepsilon^{q_7-q_2}E_{n_{22}}),\quad P={\rm diag}(P_1,P_2),$$
where $E_n$ represents the $n\times n$ identity matrix. Then the system (1.1) reads
\begin{equation}\label{(4.1)}
\left(\begin{array}{c} \dot I\\ \dot\varphi\end{array}\right)=\left(\begin{array}{c} A^0(\xi,\varepsilon)I\\ \omega^0(\xi,\varepsilon)\end{array}\right)+PG(I,\varphi,\xi,\varepsilon)
\end{equation}
with $G=\varepsilon^{q_2}{\rm col}(g_1,g_2,g_3,g_4)$.
\vskip 0.2in

a) {\bf Outline of the proof}   We are going to prove Theorem 1 by employing the KAM iteration process. By Corollary~\ref{corollary2} (see Section 3), we obtain a sequence of real analytic functions  $G^0=0, G^j=\varepsilon^{q_2}{\rm col}(g_1^j,g_2^j,g_3^j,g_4^j)(j=1,2,\cdots)$ approximating $G$ and
\begin{equation}\label{gappr}
\lim_{j\rightarrow \infty}||G^j-G||_{p,\alpha;\Omega^*\times \mathbb{T}^{n_2},\Pi}=0 \qquad {\rm for \,\, all} \,\, 0\leq p<l,
\end{equation}
\begin{equation}\label{gjappr}
|G^j-G^{j-1}|_{\mathcal{U}_{j-1},\alpha;\Pi}\leq C_0 r_{j-1}^l||G||_{l,\alpha;\Omega\times \mathbb{T}^{n_2},\Pi} \qquad {\rm for }\,j\geq 1.
\end{equation}
The definitions of $\Omega^*,\mathcal{U}_j$ etc are seen above Corollary~\ref{corollary2}. Denote $G_1^j=\varepsilon^{q_2}{\rm col}(g_1^j,g_2^j)$ and $G_2^j=\varepsilon^{q_2} {\rm col}(g_3^j,g_4^j)$. We truncate $G^1$ to its lower-degree terms
$$\mathfrak{L}(G^1)\coloneqq \left(\begin{array}{c} G_1^1(0,\varphi)+\partial_IG_1^1(0,\varphi)I\\ G_2^1(0,\varphi)\end{array}\right)\coloneqq \left(\begin{array}{c} u_0^0(\varphi)+u_1^0(\varphi)I\\ w^0(\varphi)\end{array}\right)$$
and write (4.1) as
\begin{equation}\label{(4.2)}
\left(\begin{array}{c} \dot I\\ \dot\varphi\end{array}\right)=\left(\begin{array}{c} A^0I\\ \omega^0\end{array}\right)+P\left(\begin{array}{c} u_0^0(\varphi)+u_1^0(\varphi)I+H_1^0\\ w^0(\varphi)+H_2^0\end{array}\right)+P (G-G^1),
\end{equation}
with $\partial_If(I,\varphi)$ represents the partial derivative (Jacobian matrix) of $f$ with respect to the variable $I$. Here, we drop parameters from functions and will do this also in the sequel whenever there is no confusion.

Moreover, the Cauchy inequality (see Lemma A.3 in \cite{Poschel96}) implies
\begin{equation}\label{u0est}
|u_0^0|_{r_0,\alpha;\Pi}\leq C_0 M \varepsilon^{q_2} r_0^l, \qquad |u_1^0|_{r_0,\alpha;\Pi}\leq C_0 M \varepsilon^{q_2} r_0^{l-1},
\end{equation}
\begin{equation}\label{w0est}
|w^0|_{r_0,\alpha;\Pi}\leq C_0 M \varepsilon^{q_2} r_0^l, \qquad H_1^0=O_{\mathcal{U}_1,\alpha;\Pi}(I^2),\quad H_2^0=O_{\mathcal{U}_1,\alpha;\Pi}(I)
\end{equation}
and
\begin{equation}\label{H0est}
|H_1^0|_{\mathcal{U}_1,\alpha;\Pi}\leq 2C_0 M \varepsilon^{q_2} r_0^{l-2}, \qquad |H_2^0|_{\mathcal{U}_1,\alpha;\Pi}\leq C_0 M \varepsilon^{q_2} r_0^{l-1},
\end{equation}
where $ M \varepsilon^{q_2}=||G||_{l,\alpha;\Omega\times \mathbb{T}^{n_2},\Pi}$.

We want to look for a transformation $T_1$ to eliminate the lower-degree terms of $PG^1$ such that in new coordinates the lower-degree terms of analytic part in \eqref{(4.2)} are much smaller than the old ones. Assume that at the $\nu$-th step of the process,  we have already found a coordinate transformation $T_{\nu}(\nu\geq 0$ with $T_0={\rm Id}$, the identity map) such that the system (4.1) is transformed into
$$\left(\begin{array}{c} \dot I\\ \dot\varphi\end{array}\right)=\left(\begin{array}{c} A^\nu I\\ \omega^\nu\end{array}\right)+P\left(\begin{array}{c} \tilde{u}_0^\nu(\varphi)+\tilde{u}_1^\nu(\varphi)I+\tilde{H}_1^\nu\\ \tilde{w}^\nu(\varphi)+\tilde{H}_2^\nu\end{array}\right)+P\mathfrak{D}_\nu (G\circ T_\nu-G^\nu\circ T_\nu),$$
where $\tilde{H}_1^\nu=O(I^2),\tilde{H}_2^\nu=O(I), \mathfrak{D}_\nu=P^{-1}(DT_\nu)^{-1}P$, the circle "$\circ$" indicates composition of functions and $DT$ the Jacobian matrix of $T$ with respect to coordinate variables. Then we replace $G^\nu$ with $G^{\nu+1}$ which is closer to $G$, and the above equation is rewritten as
\begin{equation}\label{eqv}
\left(\begin{array}{c} \dot I\\ \dot\varphi\end{array}\right)=\left(\begin{array}{c} A^\nu I\\ \omega^\nu\end{array}\right)+P\left(\begin{array}{c} u_0^\nu(\varphi)+u_1^\nu(\varphi)I+H_1^\nu\\ w^\nu(\varphi)+H_2^\nu\end{array}\right)+P\mathfrak{D}_\nu (G\circ T_\nu-G^{\nu+1}\circ T_\nu),
\end{equation}
where
$$\left(\begin{array}{c} u_0^\nu(\varphi)+u_1^\nu(\varphi)I\\ w^\nu(\varphi)\end{array}\right)= \left(\begin{array}{c} \tilde{u}_0^\nu(\varphi)+\tilde{u}_1^\nu(\varphi)I\\ \tilde{w}^\nu(\varphi)\end{array}\right) + \mathfrak{L}(\mathfrak{D}_\nu (G^{\nu+1}\circ T_\nu-G^{\nu}\circ T_\nu)),$$
$$H_1^\nu=O(I^2),\qquad H_2^\nu=O(I).$$
We want to construct a coordinate change $T^{\nu+1}$ to eliminate the lower-degree terms in \eqref{eqv} such that the lower-degree terms of the next step are much smaller. Repetition of this process leads to a sequence of transformation $T_\nu=T_{\nu-1}\circ T^\nu$ with $T_0={\rm Id}, \nu=1,2,\cdots$, the limit transformation of which , if converges, reduces (4.1) into a system without the lower-degree terms. Thus, we can obtain the quasi-periodic solution of (4.1). The proof of convergence is due to the following iteration lemma which describes quantitatively the KAM iteration process.
\vskip 0.2in

b) {\bf Iteration Lemma}  Before stating the iteration lemma we first introduce the iterative sequences and notations used at each iteration step. Set
$$\varepsilon_0=\varepsilon^{q_2},\qquad ||G||_{l.\alpha;\Omega\times \mathbb{T}^{n_2},\Pi}=M \varepsilon_0,$$
 $$\Omega=\{I\in \mathbb{R}^{n_1}:\, |I|<3\tilde{r}\},\quad \Omega^*=\{I\in \mathbb{R}^{n_1}:\, |I|\leq 2\tilde{r}\},\quad \Omega_0=\{I\in \mathbb{R}^{n_1}:\, |I|<\tilde{r}\}$$
with some constant $0<\tilde{r}\leq 1$. For $\nu\geq 1$, let

(i) $ r_0=\tilde{r},\, r_\nu=\tilde{r}3^{-\nu}$,
$$\mathcal{U}_\nu=\Omega^*\times \mathbb{T}^{n_2}+(3r_\nu,3r_\nu),\quad \mathcal{V}_\nu =\Omega_0\times \mathbb{T}^{n_2}+(r_\nu,r_\nu), \qquad \nu\geq 0,$$
$$ \mathcal{V}_\nu^*=\Omega_0\times \mathbb{T}^{n_2}+(2r_\nu,2r_\nu)\subset \mathcal{V}_{\nu-1}\subset \mathcal{U}_\nu;$$

(ii) $K_0=0,\quad K_\nu=[K_\nu^\prime ]+1, \quad K_\nu^\prime =3^\nu r_0^{-1}(\ln \widetilde{C}+(n_2+1)|\ln r_0|+(l+(n_2+1)\nu-\alpha)\ln3), \, \widetilde{C}=24(n_2!)n_2^{n_2}e^{-n_2}$, $[K_\nu^\prime ]$ is the integer part of $K_\nu^\prime$;

(iii) $s_0=\gamma, \quad s_\nu=\gamma (16c_1n_3\sqrt{n_2}K_\nu^{\iota+1})^{-1},\quad \Pi_\nu^{s_\nu}=\{\xi \in \mathbb{R}^{n_3}:\, {\rm dist}(\xi, \Pi_\nu)<s_\nu\};$

(iv) $\chi_\nu=r_\nu^{l-2(\alpha+1)(\iota+1)-\alpha-3},\quad X_\nu=\sum_{j=1}^{\nu}\chi_j,$\\
 the assumption $l>2(\alpha+1)(\iota+2)+\alpha\iota$ implies $X_\nu=\sum_{j=1}^{\infty}\chi_j<\frac{1}{2}$.

(v) $\delta_{\nu\mu}=\gamma^{-\mu-1}r_\nu^{l-(\alpha+\mu+2)(\iota+1)-\alpha-3}C_0M\varepsilon_0,\quad 0\leq \mu\leq \alpha;$

(vi) $f(I,\varphi,\xi)=O_{\mathcal{U},\alpha;\Pi}(I^k)$ denotes a map which is real analytic in coordinate variables $(I,\varphi)\in \mathcal{U}$, continuously differentiable up to order $\alpha$ in parameter $\xi \in \Pi$, and vanishes with $I$-derivatives up to order $k-1\geq 0$, and $f$ and its $\xi$-derivatives up to order $\alpha$ are bounded on $\mathcal{U}\times \Pi$.

\begin{lemma}\label{iterationlemma}
 (Iteration Lemma) Assume that for the equation \eqref{eqv} with $\nu\geq 0$,

{\rm (v.1)}\, (Frequency condition) let $A^\nu={\rm diag}(\varepsilon^{q_1}A_1^\nu, \varepsilon^{q_3}A_2^\nu),\quad \Lambda^\nu={\rm diag}(\varepsilon^{q_1}\Lambda_1^\nu, \varepsilon^{q_3}\Lambda_2^\nu), \quad \Lambda_1^\nu={\rm diag}(\lambda_1^\nu,\cdots,\lambda_{n_{11}}^\nu), \quad \Lambda_1^\nu={\rm diag} (\lambda_{n_{11}+1}^\nu,\cdots,\lambda_{n_{1}}^\nu)$, $A_i^\nu=B_i\Lambda_i^\nu B_i^{-1}\, (i=1,2)$ and $ \omega^\nu={\rm col}(\varepsilon^{q_5}\omega_1^\nu,\omega_2^\nu)$ satisfy, for $\varepsilon \in (0,\varepsilon^*]$,
$$\inf_{\xi\in \Pi_\nu^{s_\nu}}|\lambda_j|\geq c_0(1-X_\nu)>\frac{c_0}{2},\quad \inf_{\xi\in \Pi_\nu^{s_\nu}}|\lambda_j-\lambda_i|\geq c_0(1-X_\nu)>\frac{c_0}{2}$$
 for $i\neq j,\, 1\leq i,j\leq n_{11},\,{\rm or} \, n_{11}+1\leq i,j\leq n_1$, and
$$ ||\Lambda_i^\nu||_{\alpha;\Pi_\nu^{s_\nu}}\leq c_1(1+X_\nu)<2c_1, \qquad \parallel \partial_{\xi}^{\beta}\omega^\nu \parallel_{\Pi_\nu^{s_\nu}}\leq c_1(1+X_\nu) \varepsilon^{q_5}<2c_1\varepsilon^{q_5}, \, 1\leq |\beta|_1\leq \alpha,$$
\begin{equation}\label{freerror}
||\Lambda_i^\nu-\Lambda_i^{\nu-1}||_{\alpha;\Pi_\nu^{s_\nu}}\leq \tilde{c}_0 C_0M\varepsilon_0\varepsilon^{b_i}r_{\nu-1}^{l-(\alpha+1)(\iota+2)-1},\,
||\omega_i^\nu-\omega_i^{\nu-1}||_{\alpha;\Pi_\nu^{s_\nu}}\leq  C_0M\varepsilon_0\varepsilon^{b_{i+2}}r_{\nu-1}^{l-(\alpha+1)(\iota+2)}, \,\nu\geq 1
\end{equation}
 for $i=1,2$,  where $\tilde{c}_0$ is a positive constant, $c_0$ and $c_1$ are given in Assumption (H2), $b_1=0, b_2=q_4-q_2, b_3=q_6-q_2, b_4=q_7-q_2$;

{\rm (v.2)}\, (Small condition) the terms $u_0^\nu, u_1^\nu$ and $w^\nu$ satisfy the following estimates
$$|u_0^\nu|_{r_\nu,\alpha;\Pi_\nu^{s_\nu}}\leq 4 C_0M\varepsilon_0r_{\nu}^{l-\alpha},\quad |u_1^\nu|_{r_\nu,\alpha;\Pi_\nu^{s_\nu}}\leq  C_0M\varepsilon_0r_{\nu}^{l-(\alpha+1)(\iota+2)-1},\quad |w^\nu|_{r_\nu,\alpha;\Pi_\nu^{s_\nu}}\leq  C_0M\varepsilon_0r_{\nu}^{l-(\alpha+1)(\iota+2)},$$
 $H_1^\nu(I,\varphi,\xi)$ and $H_2^\nu(I,\varphi,\xi)$ fulfill
 \begin{equation}\label{hest}
H_1^\nu=O_{\mathcal{V}_\nu,\alpha;\Pi_\nu^{s_\nu}}(I^2),\quad H_2^\nu=O_{\mathcal{V}_\nu,\alpha;\Pi_\nu^{s_\nu}}(I),\quad |H_i^\nu-H_i^{\nu-1}|_{\mathcal{V}_\nu,\alpha;\Pi_\nu^{s_\nu}}\leq \chi_\nu C_0M\varepsilon_0
\end{equation}
for $\nu\geq 1, \,i=1,2$;

{\rm (v.3)}\, (Transformation) the transformation $T_\nu:\mathcal{V}_\nu\times\Pi_\nu^{s_\nu}\rightarrow \mathcal{U}_\nu$ is real analytic in coordinate variables $(I,\varphi)\in \mathcal{V}_\nu$ and continuously differentiable up to order $\alpha$ in the parameter $\xi \in \Pi_\nu^{s_\nu}$, satisfies
\begin{equation}\label{transest}
 |T_\nu-T_{\nu-1}|_{\mathcal{V}_\nu,\mu;\Pi_\nu^{s_\nu}}\leq (1+X_\nu)C_1 C_0M\varepsilon_0\gamma^{-\mu-1}r_{\nu}^{l-(\alpha+\mu+2)(\iota+1)-\alpha-2}<r_\nu \chi_\nu,
\end{equation}
\begin{equation}\label{dTvest}
 |P^{-1}(DT_\nu-DT_{\nu-1})P|_{\mathcal{V}_\nu,\mu;\Pi_\nu^{s_\nu}}\leq 2(1+X_\nu)C_1 C_0M\varepsilon_0\gamma^{-\mu-1}r_{\nu}^{l-(\alpha+\mu+2)(\iota+1)-\alpha-3}< \chi_\nu
\end{equation}
with $T_0={\rm Id}$ and $0\leq \mu\leq \alpha$, where $C_1$ is a constant independent of $\nu$.

Then there exists a closed set $\Pi_{\nu+1}\subset \Pi_\nu$
$$\Pi_{\nu+1}=\left\{\xi \in \Pi_\nu:\, |\sqrt{-1}\langle k,\omega^{\nu}\rangle+\langle m,\Lambda^{\nu}\rangle|\geq \gamma \varepsilon^{q_5}|k|_2^{-\iota},m\in \mathfrak{m},k\in\mathbb{Z}^{n_2}, K_{\nu}<|k|_2\leq K_{\nu+1} \right\}$$
(see Theorem 1 and (H3) for definitions of $ \mathfrak{m}$ and $\iota$, respectively) and a coordinate transformation
$$T^{\nu+1}:\mathcal{V}_{\nu+1}\times\Pi_{\nu+1}^{s_{\nu+1}}\rightarrow \mathcal{V}_{\nu+1}^*\subset \mathcal{V}_\nu \subset \mathcal{U}_{\nu+1}$$
in the form
\begin{equation}\label{vtran}
I=\rho +v_0^\nu(\phi,\xi)+v_1^\nu(\phi,\xi)\rho,\qquad \varphi=\phi+\Phi^\nu(\phi,\xi),
\end{equation}
where $\rho$ and $\phi$ are new coordinate variables, and all terms in the transformation are real analytic in $\phi$ and continuously differentiable in $\xi$ up to order $\alpha$, satisfy the estimates
\begin{equation}\label{phiv}
 |\Phi^\nu|_{2r_{\nu+1},\alpha;\Pi_{\nu+1}^{s_{\nu+1}}}\leq C_1 C_0M\varepsilon_0\gamma^{-\alpha-1}r_{\nu+1}^{l-(\alpha+1)(2\iota+3)},
\end{equation}
\begin{equation}\label{v0v}
 |v_0^\nu|_{2r_{\nu+1},\alpha;\Pi_{\nu+1}^{s_{\nu+1}}}\leq C_1 C_0M\varepsilon_0\gamma^{-\alpha-1}r_{\nu+1}^{l-(\alpha+1)(\iota+1)-\alpha},
\end{equation}
\begin{equation}\label{v1v}
 |v_1^\nu|_{2r_{\nu+1},\alpha;\Pi_{\nu+1}^{s_{\nu+1}}}\leq C_1 C_0M\varepsilon_0\gamma^{-\alpha-1}r_{\nu+1}^{l-(\alpha+1)(2\iota+3)-1}
\end{equation}
and
\begin{equation}\label{dest}
 |P^{-1}(DT^{\nu+1})^{-1}P|_{\mathcal{V}_{\nu+1},0;\Pi_{\nu+1}^{s_{\nu+1}}}<1+\chi_{\nu+1},\quad  \left|\partial_\xi^\beta \left(P^{-1}(DT^{\nu+1})^{-1}P\right)\right|_{\mathcal{V}_{\nu+1},0;\Pi_{\nu+1}^{s_{\nu+1}}}<\chi_{\nu+1}
\end{equation}
  for $ 1\leq |\beta|_1\leq \alpha$, such that the equation \eqref{eqv} is transformed into
$$\left(\begin{array}{c} \dot \rho\\ \dot\phi\end{array}\right)=\left(\begin{array}{c} A^{\nu+1} \rho\\ \omega^{\nu+1}\end{array}\right)+P\left(\begin{array}{c} u_0^{\nu+1}(\phi)+u_1^{\nu+1}(\phi)\rho+H_1^{\nu+1}\\ w^{\nu+1}(\phi)+H_2^{\nu+1}\end{array}\right)+P\mathfrak{D}_{\nu+1} (G\circ T_{\nu+1}-G^{\nu+2}\circ T_{\nu+1})$$
and the conditions (v.1)-(v.3) are satisfied by replacing $\nu$ by $\nu+1$ and $(I,\varphi)$ by $(\rho,\phi)$, respectively, where $T_{\nu+1}=T_\nu\circ T^{\nu+1}, \mathfrak{D}_{\nu+1}=P^{-1}(DT_{\nu+1})^{-1}P$.
\end{lemma}
\vskip 0.2in

c) {\bf Proof of Theorem 1}  Theorem 1 is easy to be proven by the Iteration Lemma and Inverse Approximation Lemma.

First the system (1.1) has been written in the form \eqref{(4.2)} just as \eqref{eqv} satisfying the conditions (v.1)-(v.3) with $\nu=0$ in the Iteration Lemma by Assumptions (H2) and (H3), \eqref{u0est} and \eqref{w0est}. We use the Iteration Lemma inductively to obtain a sequence of transformations $T_\nu$ mapping $\mathcal{V}_\nu\times \Pi_\nu^{s_\nu}$ into $\mathcal{V}_0$ and satisfying the estimate \eqref{transest}. Noting that $\mathcal{V}_\nu$ and $\Pi_\nu^{s_\nu}$ are exactly regarded as those neighbourhoods of the open convex set $\Omega_0\times\mathbb{T}^{n_2}\subset \mathbb{R}^{n_1+n_1}$ and closed subset $\Pi_\gamma\subset \Pi$, respectively, and $r_\nu^{l_1}/s_\nu^{\mu}\rightarrow 0$ as $\nu\rightarrow \infty$ $(l_1=l-(\alpha+\mu+2)(\iota+1)-\alpha-2$ and the positive integer $\mu\leq\alpha)$  by the definition of $s_\nu$, the Inverse Approximation Lemma and Condition (v.3) imply that for every $\xi \in \Pi_\gamma$, the limit map $T=\lim_{\nu\rightarrow \infty}T_\nu$ exists in $C^{p,\mu-1}(\Omega_0\times\mathbb{T}^{n_2},\Pi_\gamma)$ for $0\leq p<l_1$ and $T:\Omega_0\times\mathbb{T}^{n_2}\times\Pi_\gamma\rightarrow \Omega^*\times\mathbb{T}^{n_2}$ for sufficiently small $\varepsilon$, and is of the form
$$T:\qquad I=\rho +V_0(\phi,\xi)+V_1(\phi,\xi)\rho,\quad \varphi=\phi+\Phi(\phi,\xi)$$
by \eqref{vtran}, which is of class $C^\mu$ in $\xi\in\Pi_{\gamma}$ in the sense of whitney and of class $\hat{C}^{l_1}$ in $\varphi\in \mathbb{T}^{n_2}$ together with derivatives up to order $\mu-1$ with respect to $\xi$ for $0< \mu\leq \alpha$. Moreover, by \eqref{gappr}, we obtain
\begin{equation}\label{remest}
\lim_{\nu\rightarrow \infty}||G\circ T_\nu-G^{\nu+1}\circ T_\nu||_{p,\mu-1;\Omega_0\times \mathbb{T}^{n_2},\Pi_\gamma}=0 \qquad {\rm for \,\, all} \,\, 0\leq p<l_1,\,0<\mu\leq \alpha
\end{equation}
and by Condition (v.3) and \eqref{dest},
\begin{equation}\label{ddest}
 |\mathfrak{D}_\nu-\mathfrak{D}_{\nu-1}|_{\mathcal{V}_\nu,\mu;\Pi_\nu^{s_\nu}}\leq C_2 C_0M\varepsilon_0\gamma^{-\mu-1}r_{\nu}^{l-(\alpha+\mu+2)(\iota+1)-\alpha-3},
\end{equation}
where $C_2$ is a constant independent of $\nu,\gamma$ and $\varepsilon_0$. It follows from \eqref{hest}, \eqref{remest} and \eqref{ddest} that System (4.1) is transformed by $T$ into the system
\begin{equation}\label{infeq}
 \left\{\begin{array}{rl}
\dot{\rho} & =A^*(\xi)\rho+P_1O(\rho^2)\\
\dot{\phi} & = \omega^*(\xi)+P_2O(\rho)
\end{array}\right.
 \end{equation}
for $(\rho,\phi)\in \Omega_0\times\mathbb{T}^{n_2}, \xi\in\Pi_\gamma$, where $A^*(\xi)={\rm diag}(\varepsilon^{q_1}A_1^*(\xi),\varepsilon^{q_3}A_2^*(\xi)), \omega^*(\xi)={\rm col}(\varepsilon^{q_5}\omega_1^*(\xi),\omega_2^*(\xi))$, $A_i^*=\lim_{\nu\rightarrow \infty}A_i^\nu$ and $\omega_i^*=\lim_{\nu\rightarrow \infty}\omega_i^\nu (i=1,2)$ exist by Condition (v.1) and are of class $C^\alpha$ in $\xi\in\Pi_\gamma$ in the sense of Whitney by the Inverse Approximation Lemma since $r_\nu^{l-(\alpha+1)(\iota+2)-1}s_\nu^{-\alpha}\rightarrow 0$ as $\nu\rightarrow \infty$. Thus, we obtain the quasi-periodic invariant torus of (4.1)
$$ I=V_0(\phi,\xi),\quad \varphi=\phi+\Phi(\phi,\xi),\quad \phi=\omega^*(\xi) t+\phi_0$$
satisfying the estimates \eqref{torusest}, \eqref{freqest} and \eqref{frees} by Conditions (v.1) and (v.3). The rest of Theorem 1 can be derived immediately from the Iteration Lemma. \hskip 0.4in $\blacksquare$
\vskip 0.2in

d) {\bf Proof of Iteration Lemma} To simplify the notation, we denote quantities referring to $\nu+1$ with $+$ such as $u^{\nu+1}$ by $u^+$, $r_{\nu+1}$ by $r_+$, and those referring to $\nu$ without the $\nu$ such as $u^\nu$ by $u$, $r_\nu$ by $r$. Substituting the transformation $T^+$ into \eqref{eqv}, the transformation $T^+$ will be obtained by solving the homological equations
\begin{eqnarray}
& &\partial_{\phi}v_0\cdot \omega-A v_0 =P_1\Gamma_{K_+} u_0(\phi),\label{v0eq}\\
& &\partial_{\phi}v_1\cdot \omega+v_1A-Av_1 =P_1(\Gamma_{K_+}u_1(\phi)-B{\rm diag}(B^{-1}\widehat{u_1}(0)B)B^{-1}),\label{v1eq}\\
& &\partial_{\phi}\Phi \cdot \omega =P_2(\Gamma_{K_+} w(\phi)-\widehat{w}(0)),\label{phieq},
\end{eqnarray}
where $B={\rm diag}(B_1,B_2), {\rm diag}(B^{-1}\widehat{u_1}(0)B)$ denotes a diagonal matrix whose elements are the diagonal elements of $B^{-1}\widehat{u_1}(0)B$, $\widehat{u_1}(0)$ and $\widehat{w}(0)$ denote the mean values (that is, the zero-order coefficients of the Fourier series expansions) of $u_1$ and $w$ over $ \mathbb{T}^{n_2}$, respectively, $\Gamma_{K_+}$ is the truncation operator of the Fourier series expansions defined in Lemma~\ref{lemmaA5} and the notation $\partial_\phi f\cdot \omega$ means $\partial_\phi f\cdot \omega=\sum_{j=1}^{n_2}t_j\frac{\partial f}{\partial \phi_j}$ for $\omega={\rm col}(t_1,\cdots,t_{n_2})$. Here, the homological equations are approximated by truncating the Fourier series expansions of $u_0,u_1$ and $w$ so that the solutions are defined on an open set of parameters. This idea is due to Arnol'd\cite{Arn63a} and P\"{o}schel\cite{Pos82}.

d1) Solutions of \eqref{v0eq}-\eqref{phieq} and estimates. Set
$$\Pi_+=\left\{\xi \in \Pi_\nu:\, |\sqrt{-1}\langle k,\omega\rangle+\langle m,\Lambda\rangle|\geq \gamma \varepsilon^{q_5}|k|_2^{-\iota},m\in \mathfrak{m}, K<|k|_2\leq K_{+} \right\}$$
and
$$\Pi_+^{s_+}=\{\xi\in \mathbb{R}^{n_3}:\, {\rm dist}(\xi,\Pi_+)<s_+\}\subset \Pi_\nu^{s_\nu}.$$

\begin{lemma}\label{lemma4.1}
For every $\xi\in\Pi_+^{s_+}$, we have
\begin{equation}\label{smalld1}
|\sqrt{-1}\langle k,\omega(\xi)\rangle+\langle m,\Lambda(\xi)\rangle|\geq \frac{1}{4}\gamma \varepsilon^{q_5}|k|_2^{-\iota},\quad 0<|k|_2\leq K_{+}, m\in \mathfrak{m}.
\end{equation}
\end{lemma}

{\bf Proof} We first prove
\begin{equation}\label{smalld2}
|\sqrt{-1}\langle k,\omega(\xi)\rangle+\langle m,\Lambda(\xi)\rangle|\geq \frac{1}{2}\gamma \varepsilon^{q_5}|k|_2^{-\iota},\quad 0<|k|_2\leq K_{+}, m\in \mathfrak{m}
\end{equation}
for every $\xi\in\Pi_+$. Noting the fact that $K_j^{\iota+1}r_j^{l-(\alpha+1)(\iota+2)-2} \rightarrow 0$ as $j\rightarrow \infty$, the \eqref{freerror} implies that for $0<|k|_2\leq K_j,1\leq j \leq \nu$,
\begin{eqnarray*}
|\sqrt{-1}\langle k,\omega^j(\xi)-\omega^{j-1}(\xi)\rangle+\langle m,\Lambda^j(\xi)-\Lambda^{j-1}(\xi)\rangle| & \leq & \varepsilon^{q_5}(\sqrt{n_2}K_j+2\tilde{c}_0) C_0M\varepsilon_0r_{j-1}^{l-(\alpha+1)(\iota+2)-1}\\
& < & r_j\gamma \varepsilon^{q_5}K_j^{-\iota}
\end{eqnarray*}
for sufficiently small $\varepsilon_0$. As for $K_{j-1}<|k|_2\leq K_j$, $|\sqrt{-1}\langle k,\omega^{j-1}(\xi)\rangle+\langle m,\Lambda^{j-1}(\xi)\rangle|\geq \gamma \varepsilon^{q_5}|k|_2^{-\iota}$, hence,
\begin{eqnarray*}
& & |\sqrt{-1}\langle k,\omega^\nu(\xi)\rangle+\langle m,\Lambda^\nu(\xi)\rangle|\\
& & \quad \geq |\sqrt{-1}\langle k,\omega^{j-1}(\xi)\rangle+\langle m,\Lambda^{j-1}(\xi)\rangle|-
\sum_{i=j}^{\nu}|\sqrt{-1}\langle k,\omega^i(\xi)-\omega^{i-1}(\xi)\rangle+\langle m,\Lambda^i(\xi)-\Lambda^{i-1}(\xi)\rangle| \\
& & \quad \geq \gamma \varepsilon^{q_5}(|k|_2^{-\iota}-\sum_{i=j}^{\nu}r_iK_i^{-\iota})\geq \frac{1}{2}\gamma \varepsilon^{q_5}|k|_2^{-\iota},
\end{eqnarray*}
which implies \eqref{smalld2}.

For every $\xi\in\Pi_+^{s_+}\subset \Pi_\nu^{s_\nu}$, there is $\xi_0\in \Pi_+$ such that $|\xi-\xi_0|<s_+$. The condition (v.1) and \eqref{smalld2} imply
\begin{eqnarray*}
& & |\sqrt{-1}\langle k,\omega(\xi)\rangle+\langle m,\Lambda(\xi)\rangle|\\
& & \quad \geq |\sqrt{-1}\langle k,\omega(\xi_0)\rangle+\langle m,\Lambda(\xi_0)\rangle|-
|\sqrt{-1}\langle k,\omega(\xi)-\omega(\xi_0)\rangle+\langle m,\Lambda(\xi)-\Lambda(\xi_0)\rangle| \\
& & \quad \geq \frac{1}{2}\gamma \varepsilon^{q_5}|k|_2^{-\iota}-2c_1n_3\varepsilon^{q_5}(\sqrt{n_2}|k|_2+2)s_+\\
& & \quad \geq \frac{1}{4}\gamma \varepsilon^{q_5}|k|_2^{-\iota}. \hskip 0.4in \blacksquare
\end{eqnarray*}
\vskip 0.2in

 The procedure of solving \eqref{v0eq}-\eqref{phieq} is standard in KAM theory. Expanding the functions into the Fourier series in $\phi$, and substituting in \eqref{v0eq}-\eqref{phieq} and comparing coefficients of the term $e^{\sqrt{-1}\langle k,\phi\rangle}$, one obtain the solutions
\begin{eqnarray}
& &v_0(\phi)=P_1\sum_{|k|_2\leq K_+} B (\sqrt{-1}\langle k,\omega\rangle-\Lambda)^{-1}B^{-1}\widehat{u_0}(k)e^{\sqrt{-1}\langle k,\phi\rangle},\label{v0}\\
& &v_1(\phi)=P_1\sum_{|k|_2\leq K_+} BV_1(k)B^{-1}e^{\sqrt{-1}\langle k,\phi\rangle},\label{v1}\\
& &\Phi(\phi)=P_2\sum_{0<|k|_2\leq K_+} (\sqrt{-1}\langle k,\omega\rangle)^{-1}\widehat{w}(k)e^{\sqrt{-1}\langle k,\phi\rangle},\label{phi}
\end{eqnarray}
where
$$(V_1(k))_{ij}=\left\{\begin{array}{ll} (\sqrt{-1}\langle k,\omega\rangle+\varepsilon^{a_j}\lambda_j-\varepsilon^{a_i}\lambda_i)^{-1}(\widehat{U_1}(k))_{ij},\,& |k|+|i-j|\neq 0\\
0, & |k|+|i-j|=0,
\end{array}\right.$$
$U_1(\phi)=B^{-1}u_1(\phi)B,\, a_i=q_1$ if $1\leq i\leq n_{11},\, =q_3$ if $n_{11}+1\leq i\leq n_1$, $(V_1(k))_{ij}$ and $(\widehat{U_1}(k))_{ij}$ represent elements of the matrices $V_1(k)$ and $\widehat{U_1}(k)$, respectively, and $\hat{u}(k)$ is the $k$-order coefficients of the Fourier series expansions of $u$. Hence, Lemma~\ref{lemma4.1}, \eqref{v0}-\eqref{phi} and Conditions (v.1)-(v.2) imply that $v_0,v_1$ and $\Phi$ are real analytic in $\phi \in \mathcal{W}\coloneqq \mathbb{T}^{n_2}+2 r_+$, continuously differentiable up to order $\alpha$ in $\xi\in\Pi_+^{s_+}$. Meanwhile using Lemma~\ref{lemmaA7}, one easily gets the estimates \eqref{phiv}-\eqref{v1v} (we denote $|\cdot|_{r,\mu; \Pi_+^{s_+}}$ by  $|\cdot|_{r,\mu; s_+}$, $|\cdot|_{r,\mu; \Pi_\nu^{s_\nu}}$ by  $|\cdot|_{r,\mu; s}$ for simplification) and
\begin{equation}\label{dv0est}
 |v_0|_{2r_+,\mu;s_+}\leq C_1\delta_{+\mu}r_+^{(\alpha+1)(\iota+1)+3}, \quad |P_1^{-1}\partial_\phi v_0|_{2r_+,\mu;s_+}\leq C_1\varepsilon^{-q_5}\delta_{+\mu}r_+^{(\alpha+1)(\iota+1)+2},
 \end{equation}
\begin{equation}\label{dv1est}
|v_1|_{2r_+,\mu; s_+} \leq  C_1\delta_{+\mu}r_+, \quad |P_1^{-1}v_1P_1|_{2r_+,\mu;s_+}\leq C_1\delta_{+\mu}r_+,
\quad |P_1^{-1}\partial_\phi v_1|_{2r_+,\mu;s_+}\leq C_1\varepsilon^{-q_5}\delta_{+\mu}
\end{equation}
and
\begin{equation}\label{dphiest}
 |\Phi|_{2r_+,\mu;s_+}\leq C_1\delta_{+\mu}r_+^2, \quad |P_2^{-1}\partial_\phi \Phi|_{2r_+,\mu;s_+}\leq C_1\varepsilon^{-q_5}\delta_{+\mu}r_+
\end{equation}
for $0\leq \mu\leq \alpha$ and an appropriate choice of the constant $C_1$ independent of $\nu$.

It is easy to see that when the $\varepsilon_0$ is sufficiently small, the transformation $T^+$ maps $\mathcal{V}_+$ into $\mathcal{V}_+^*\subset \mathcal{V}$  and $\mathcal{V}_+^*$ into $\mathcal{V}$, respectively, and
\begin{equation}\label{vTest}
|T^+-{\rm Id}|_{\mathcal{V}_+^*,0;s_+}\leq C_1 r_+\delta_{+0},\quad |\partial_\xi^\beta T^+|_{\mathcal{V}_+^*,0;s_+}\leq C_1 r_+\delta_{+\mu},\quad 1\leq |\beta|_1=\mu \leq \alpha.
\end{equation}

d2) Proof of \eqref{dest}.  Corresponding to the transformation $T^+$, we have its Jacobian matrix
\begin{equation}\label{dT}
DT^+=\left(\begin{array}{cc} E_{n_1}+v_1 & \partial_\phi v_0+\partial_\phi v_1\rho\\
0 & E_{n_2}+\partial_\phi \Phi\end{array}\right)
\end{equation}
and the inverse
\begin{equation}\label{dTin}
\left(DT^+\right)^{-1}=\left(\begin{array}{cc} (E_{n_1}+v_1)^{-1} & -(E_{n_1}+v_1)^{-1}(\partial_\phi v_0+\partial_\phi v_1\rho)(E_{n_2}+\partial_\phi \Phi)^{-1}\\
0 & (E_{n_2}+\partial_\phi \Phi)^{-1}\end{array}\right).
\end{equation}
Thus, \eqref{dv0est}-\eqref{dphiest} imply
\begin{equation}\label{vdTest}
|P^{-1}DT^+P|_{\mathcal{V}_+,0;s_+}\leq 1+\chi_+ r_+^{\alpha(\iota+1)},\quad |P^{-1}(DT^+-E)P|_{\mathcal{V}_+,\mu;s_+}\leq C_1\delta_{+\mu}\quad {\rm for}\,1\leq \mu\leq \alpha.
\end{equation}
Noting that the derivatives of $DT^+$ with respect to the parameter $\xi$ is sufficiently small and that for a matrix $M(\xi)$ with a small norm, differentiating the left- and right-hand sides of $(E+M(\xi))^{-1}(E+M(\xi))=E$ and using the Leibniz formula, we find
$$\partial_\xi^\beta (E+M(\xi))^{-1}=-\sum_{k<\beta}\left(\begin{array}{c} \beta\\k\end{array}\right)\partial_\xi^k (E+M(\xi))^{-1}\cdot \partial_\xi^{\beta-k} (E+M(\xi))\cdot (E+M(\xi))^{-1},$$
where $k,\beta\in \mathbb{Z}_+^{n_3}, \left(\begin{array}{c} \beta\\k\end{array}\right)=\left(\begin{array}{c} \beta_1\\k_1\end{array}\right)\cdots \left(\begin{array}{c} \beta_{n_3}\\k_{n_3}\end{array}\right)$, $E$ is the identity matrix, the estimates \eqref{dv1est} and \eqref{dphiest} imply
\begin{eqnarray}
& & |P_1^{-1}(E_{n_1}+v_1)^{-1}P_1|_{2r_+,0;s_+}\leq 1+\chi_+ r_+^{\alpha(\iota+1)+1},  \label{p1v1in}\\
& & |P_2^{-1}(E_{n_2}+\partial_\phi \Phi)^{-1}P_2|_{2r_+,0;s_+}\leq 1+\chi_+ r_+^{\alpha(\iota+1)+1},  \label{p2phiin}\\
& & |\partial_\xi ^\beta (P_1^{-1}(E_{n_1}+v_1)^{-1}P_1)|_{2r_+,0;s_+}\leq 2(1+r_+\chi_+)^2 |P_1^{-1}v_1P_1|_{2r_+,|\beta|_1;s_+},  \label{dp1v1in}\\
& & |\partial_\xi ^\beta (P_2^{-1}(E_{n_2}+\partial_\phi \Phi)^{-1}P_2)|_{2r_+,0;s_+}\leq 2(1+r_+\chi_+)^2 |P_2^{-1}\partial_\phi \Phi P_2|_{2r_+,|\beta|_1;s_+}  \label{dp2phiin}
\end{eqnarray}
for $1\leq |\beta|_1\leq \alpha$ and sufficiently small $\varepsilon_0$. Hence, \eqref{dest} follows from \eqref{dTin}, \eqref{dv0est}, \eqref{dv1est} and \eqref{p1v1in}-\eqref{dp2phiin}. Moreover, we have
\begin{equation}\label{ddest1}
 \left.\begin{array}{r}|P_1^{-1}(E_{n_1}+v_1)^{-1}P_1F|_{\mathcal{V}_+,\alpha;s_+} \\|P_2^{-1}(E_{n_2}+\partial_\phi \Phi)^{-1}P_2F|_{\mathcal{V}_+,\alpha;s_+} \\ |P^{-1}(DT^+)^{-1}PF|_{\mathcal{V}_+,\alpha;s_+} \end{array}\right\} <(1+r_+\chi_+)|F|_{\mathcal{V}_+,\alpha;s_+}
 \end{equation}
 for suitable $F$ which is real analytic in $\mathcal{V}_+$ and continuously differentiable up to order $\alpha$ in $\xi \in \Pi_+^{s_+}$.

 d3) We proceed to verify (v.3) with $\nu+1$ replacing $\nu$. As the transformation $T^+:\, \mathcal{V}_+\times \Pi_+^{s_+}\rightarrow \mathcal{V}_+^*$ (or $ \mathcal{V}_+^*\times \Pi_+^{s_+}\rightarrow \mathcal{V}$) is real analytic in coordinate variables and continuously differentiable up to order $\alpha$ in parameter $\xi$, so is $T_+=T\circ T^+$.

 We first prove \eqref{transest} and \eqref{dTvest} inductively. For $\nu=1$, by \eqref{dv0est}-\eqref{dT}, it implies
$$ |T_1-{\rm Id}|_{\mathcal{V}_1,\mu;s_1}\leq C_1 r_1\delta_{1\mu}<r_1\chi_1,\quad |DT_1-E|_{\mathcal{V}_1,\mu;s_1}\leq C_1\delta_{1\mu}<\chi_1$$
and
$$|P^{-1}(DT_1-E)P|_{\mathcal{V}_1,\mu;s_1}\leq C_1\delta_{1\mu}<\chi_1$$
for $0\leq \mu\leq\alpha$. Assume that at the $\nu$-th step we have
$$ |T-T_{\nu-1}|_{\mathcal{V},\mu;s}\leq (1+X)C_1 r\delta_{\nu\mu}<r\chi,\quad |DT-DT_{\nu-1}|_{\mathcal{V},\mu;s}\leq(1+X) C_1\delta_{\nu\mu}<\chi$$
and
$$ |P^{-1}(DT-DT_{\nu-1})P|_{\mathcal{V},\mu;s}\leq 2(1+X) C_1\delta_{\nu\mu}<\chi $$
for $0\leq \mu\leq\alpha$, here we have omitted the subscript $\nu$ from the quantities referring to $\nu$. Then in view of the induction assumptions we obtain
\begin{equation}\label{dTest}
|DT|_{\mathcal{V},0;s}\leq 1+X,\quad |\partial_\xi^\beta DT|_{\mathcal{V},0;s}\leq X \quad {\rm for}\,1\leq |\beta|_1\leq \alpha.
\end{equation}
Combining \eqref{vTest},\eqref{dTest} and Lemma~\ref{lemmaA6} (i) we get
\begin{equation}\label{T+est}
|T_+-T|_{\mathcal{V}_+^*,\mu;s_+} =|T\circ T^+-T|_{\mathcal{V}_+^*,\mu;s_+} \leq (1+X_+)C_1 r_+\delta_{+\mu}<r_+\chi_+ \quad {\rm for}\,0\leq \mu \leq \alpha,
\end{equation}
which, together with the Cauchy inequality, implies
\begin{equation}\label{dT+est}
|DT_+-DT|_{\mathcal{V}_+,\mu;s_+}  \leq (1+X_+)C_1 \delta_{+\mu}<\chi_+ \quad {\rm for}\,0\leq \mu \leq \alpha.
\end{equation}
Similarly, we have
\begin{equation}\label{dpTest}
|P^{-1}DT P|_{\mathcal{V},0;s} \leq 1+X,\quad |\partial_\xi^\beta (P^{-1}DT P)|_{\mathcal{V},0;s} \leq X \quad {\rm for}\,1\leq |\beta|_1 \leq \alpha
\end{equation}
and
\begin{equation}\label{ddpTest}
|(P^{-1}DT P)\circ T^+-P^{-1}DTP|_{\mathcal{V}_+,\mu;s_+}  \leq (1+\chi_+)(1+X)C_1 \delta_{+\mu} \quad {\rm for}\,0\leq \mu \leq \alpha
\end{equation}
by \eqref{vTest},\eqref{dpTest} and Lemma~\ref{lemmaA6} (i).  Based on the observation
$$P^{-1}(DT_+-DT)P=((P^{-1}DT P)\circ T^+-P^{-1}DTP)(P^{-1}DT^+ P)+P^{-1}DT P(P^{-1}(DT^+-E) P),$$
from \eqref{ddpTest}, \eqref{vdTest}, \eqref{dpTest} and the Leibniz formula, it follows
$$|P^{-1}(DT_+-DT)P|_{\mathcal{V}_+,\mu;s_+}  \leq 2(1+X_+)C_1 \delta_{+\mu}<\chi_+$$
 for $0\leq \mu \leq \alpha$ and sufficiently small $\varepsilon_0$. Thus, we have proved \eqref{transest} and \eqref{dTvest} with $\nu+1$.

 Now, we show that $T_+$ maps $\mathcal{V}_+$ into $\mathcal{U}_+$. Noting the expression of $T_+$ in angle variable direction is independent of $\rho$, we set $T_+(\rho,\phi)={\rm col}(v_+(\rho,\phi),\Phi_+(\phi))$. In view of the induction hypotheses, \eqref{T+est} implies
 $$|T_+-{\rm Id}|_{\mathcal{V}_+^*,0;s_+}\leq \sum_{j=1}^{\nu+1}|T_j-T_{j-1}|_{\mathcal{V}_+^*,0;s_+}\leq \sum_{j=1}^{\nu+1}r_j<r_0,$$
 and \eqref{dT+est} implies
\begin{equation}\label{dT+est1}
|DT_+|_{\mathcal{V}_+,0;s_+}  < 1+X_+<2.
\end{equation}
Hence, the first component of $T_+$ is mapped into $\Omega^*+r_+$. For $\phi$ with $|{\rm Im} \phi|<r_+$, there exists a $\phi_0\in \mathbb{T}^{n_2}$ such that $|\phi-\phi_0|<r_+$. Therefore, the $\Phi_+$ being real analytic and \eqref{dT+est1} imply
$$|{\rm Im} \Phi_+(\phi)|=|{\rm Im} (\Phi_+(\phi)- \Phi_+(\phi_0))|\leq |DT_+|_{\mathcal{V}_+,0;s_+} |\phi-\phi_0|<2r_+.$$
Thus, for $\xi \in \Pi_+^{s_+}$,  $T_+$ maps $\mathcal{V}_+$ into $\Omega^*\times \mathbb{T}^{n_2}+(2r_+,2r_+)\subset \mathcal{U}_+$, as claimed.
Furthermore, let $\mathfrak{G}(\rho,\phi)=\mathfrak{D}_+(G^{\nu+2}\circ T_+-G^{\nu+1}\circ T_+)$, then \eqref{gjappr}, \eqref{transest}, \eqref{dTvest}, Lemma~\ref{lemmaA6} (ii) and the Cauchy inequality imply
\begin{equation}\label{G+est}
|\mathfrak{G}|_{\mathcal{V}_+,\alpha;s_+} \leq 2r_+^{l-\alpha}C_0M\varepsilon_0
\end{equation}
and
\begin{equation}\label{dG+est}
|\partial_\rho\mathfrak{G}|_{\mathcal{V}_+,\alpha;s_+} \leq 2r_+^{l-\alpha-1}C_0M\varepsilon_0
\end{equation}
for sufficiently small $\varepsilon_0$.

d4) Estimates of remainder terms. Denote
$$ W(I,\varphi)={\rm col}(u_0(\varphi)+u_1(\varphi)I,w(\varphi)), \quad A^+=A+P_1\widetilde{A}, \quad \Lambda^+=\Lambda+P_1\widetilde{\Lambda},\quad \omega^+=\omega+P_2\widetilde{\omega},$$
where
$$\widetilde{A}=B({\rm diag}(B^{-1}\widehat{u_1}(0)B))B^{-1},\qquad \widetilde{\Lambda}={\rm diag}(B^{-1}\widehat{u_1}(0)B),\qquad \widetilde{\omega}=\widehat{w}(0).$$
Then the assumption (H2) implies that there is a constant $\tilde{c}_0\geq 1$ such that
\begin{equation}\label{A+est}
|{\rm diag}(B^{-1}\widehat{u_1}(0)B)|_{\alpha;s_+} \leq \tilde{c}_0 |u_1|_{r,\alpha;s},\qquad   |\widetilde{A} |_{\alpha;s_+} \leq \tilde{c}_0 |u_1|_{r,\alpha;s}.
\end{equation}
By Assumption (H1) and Condition (v.2), it is easy to see $A^+, \Lambda^+$ and $\omega^+$ satisfy (v.1) with $\nu$ replaced by $\nu+1$. We have found the transformation $T^+$ which transforms the equation \eqref{eqv}, by using \eqref{v0eq}-\eqref{phieq}, into the following one in the new variables
\begin{eqnarray*}
\left(\begin{array}{c} \dot \rho\\ \dot\phi\end{array}\right)& = & \left(\begin{array}{c} A^+ \rho\\ \omega^+\end{array}\right)+P \mathfrak{D}^+\left\{ P^{-1}(E-DT^+)P \left(\begin{array}{c} \tilde{A}\rho\\ \tilde{\omega}\end{array}\right)+({\rm Id}-\Gamma_{K_+})W(\rho,\phi)\right.\\
& & \left.+W\circ T^+(\rho,\phi)- W(\rho,\phi)+\left(\begin{array}{c} H_1\\ H_2\end{array}\right)\circ T^+ (\rho,\phi)\right\}+P\mathfrak{G}(\rho,\phi) +P\mathfrak{D}_+ (G-G^{\nu+2})\circ T_+(\rho,\phi),
\end{eqnarray*}
where $\mathfrak{D}^+=(P^{-1}DT^+P)^{-1}, \mathfrak{D}_+=(P^{-1}DT_+P)^{-1}$. We use the notation $Lf(\rho,\phi)$ to denote the linear part of a function $f$ in $\rho$, that is
$$Lf(\rho,\phi)=f(0,\phi)+\partial_\rho f(0,\phi)\rho$$
and denote $\mathfrak{G}(\rho,\phi)={\rm col}(\mathfrak{G}_1(\rho,\phi),\mathfrak{G}_2(\rho,\phi))$, rewrite the above equation in the form of \eqref{eqv}
\begin{equation}\label{eqv+}
\left(\begin{array}{c} \dot \rho\\ \dot\phi\end{array}\right)=\left(\begin{array}{c} A^+ \rho\\ \omega^+\end{array}\right)+P\left(\begin{array}{c} u_0^+(\phi)+u_1^+(\phi)\rho+H_1^+\\ w^+(\phi)+H_2^+\end{array}\right)+P\mathfrak{D}_+ (G-G^{\nu+2})\circ T_+(\rho,\phi),
\end{equation}
where
\begin{eqnarray}
& & w^+  =  P_2^{-1}(E_{n_2}+\partial_\phi\Phi)^{-1}P_2\left[-P_2^{-1}\partial_\phi\Phi P_2 \tilde{\omega}+({\rm Id}-\Gamma_{K_+})w+w(\phi+\Phi)-w(\phi)\right. \nonumber\\
& & \qquad \left.+H_2(v_0,\phi+\Phi)+\mathfrak{G}_2(0,\phi)\right],\label{w+exp}\\
& & u_0^+  =  P_1^{-1}(E_{n_1}+v_1)^{-1}P_1\left[-P_1^{-1}\partial_\phi v_0 P_2( \tilde{\omega}+w^+(\phi))+({\rm Id}-\Gamma_{K_+})u_0+u_0(\phi+\Phi)-u_0(\phi)\right. \nonumber\\
& & \qquad \left.+u_1(\phi+\Phi)v_0+H_1(v_0,\phi+\Phi)+\mathfrak{G}_1(0,\phi)\right],\label{u0+exp}\\
& & u_1^+  =  P_1^{-1}(E_{n_1}+v_1)^{-1}P_1\left[-P_1^{-1}v_1P_1\widetilde{A}-P_1^{-1}\partial_\phi v_1 P_2( \tilde{\omega}+w^+(\phi))+({\rm Id}-\Gamma_{K_+})u_1\right. \nonumber\\
 & & \qquad \left.+u_1(\phi+\Phi)-u_1(\phi)+u_1(\phi+\Phi)v_1+\partial_I H_1(v_0,\phi+\Phi)(E_{n_1}+v_1)\right. \nonumber\\
& & \qquad \left.+\partial_\rho\mathfrak{G}_1(0,\phi)-P_1^{-1}\partial_\phi v_0(E_{n_2}+\partial_\phi\Phi)^{-1} P_2Q \right],\label{u1+exp}\\
& & Q = \partial_I H_2(v_0,\phi+\Phi)(E_{n_1}+v_1)+\partial_\rho\mathfrak{G}_2(0,\phi),\nonumber
\end{eqnarray}
$$\left(\begin{array}{c} H_1^+\\ H_2^+\end{array}\right)  = \mathfrak{D}^+({\rm Id}-L)\left(\begin{array}{c} H_1\circ T^+ +\mathfrak{G}_1\\ H_2\circ T^+ +\mathfrak{G}_2 \end{array}\right)+
\left(\begin{array}{c}-P_1^{-1}(E_{n_1}+v_1)^{-1}(\partial_\phi v_1\rho) (E_{n_2}+\partial_\phi\Phi)^{-1} P_2Q\rho \\P_2^{-1}(E_{n_2}+\partial_\phi\Phi)^{-1} P_2Q\rho   \end{array}\right),$$
or in another form,
\begin{eqnarray}
& & H_2^+  =  P_2^{-1}(E_{n_2}+\partial_\phi\Phi)^{-1}P_2\left[H_2\circ T^+(\rho,\phi)+ \mathfrak{G}_2(\rho,\phi)-H_2\circ T^+(0,\phi)-\mathfrak{G}_2(0,\phi)\right],\quad \label{H2+exp}\\
& & H_1^+  =  P_1^{-1}(E_{n_1}+v_1)^{-1}P_1\left[({\rm Id}-L) (H_1\circ T^+(\rho,\phi) +\mathfrak{G}_1(\rho,\phi))\right. \nonumber \\
& & \qquad \left.+P_1^{-1}\partial_\phi v_0(E_{n_2}+\partial_\phi\Phi)^{-1} P_2Q\rho -P_1^{-1}(\partial_\phi v_0+\partial_\phi v_1\rho) P_2H_2^+\right].\label{H1+exp}
\end{eqnarray}
By \eqref{H0est} and \eqref{hest}, we have
\begin{equation}\label{Hiest}
|H_i|_{\mathcal{V},\alpha;s} \leq XC_0M\varepsilon_0+|H_i^0|_{\mathcal{U}_1,\alpha;\Pi}< \widetilde{C}_0C_0M\varepsilon_0,\quad i=1,2,
\end{equation}
where $\widetilde{C}_0=2+\sum_{j\geq 1} \chi_j<\infty.$

Now, we proceed to prove (v.2) for $\nu+1$ and first estimate the three terms $u_0^+,u_1^+$ and $w^+$. We will use $C_\alpha$ to denote a constant only depending on $\alpha$. By using the Taylor expansions of $H_1$ and $H_2$, Lemma~\ref{lemmaA6} (ii), \eqref{dv0est}-\eqref{dphiest}, \eqref{hest}, \eqref{Hiest}, \eqref{G+est} and the Cauchy inequality we find
\begin{eqnarray}
 |H_1(v_0,\phi+\Phi)|_{r_+,\alpha;s_+} &\leq & C_\alpha \widetilde{C}_0C_0M\varepsilon_0 r_+^{-2}|v_0|^2_{r_+,\alpha;s_+}\ll C_0M\varepsilon_0 r_+^{l-\alpha},\label{u0+est1}\\
 |H_2(v_0,\phi+\Phi)|_{r_+,\alpha;s_+}&\leq & C_\alpha \widetilde{C}_0C_0M\varepsilon_0 r_+^{-1}|v_0|_{r_+,\alpha;s_+}\ll C_0M\varepsilon_0 r_+^{l-(\alpha+1)(\iota+2)},\label{w+est1}\\
|Q|_{r_+,\alpha;s_+} & \leq & C_\alpha \widetilde{C}_0C_0M\varepsilon_0 r_+^{-1} \label{Qest}
\end{eqnarray}
and
\begin{equation}\label{u1+est1}
|\partial_I H_1(v_0,\phi+\Phi)(E_{n_1}+v_1) |_{r_+,\alpha;s_+}\leq  C_\alpha \widetilde{C}_0C_0M\varepsilon_0 r_+^{-2}|v_0|_{r_+,\alpha;s_+}\ll C_0M\varepsilon_0 r_+^{l-(\alpha+1)(\iota+2)-1}
\end{equation}
for sufficiently small $\varepsilon_0$. From the Cauchy inequality, Condition (v.2), Lemma~\ref{lemmaA6} (i) and \eqref{dphiest}, it follows
\begin{eqnarray}
& & |u_0(\phi+\Phi)-u_0(\phi)|_{r_+,\alpha;s_+}\leq C_\alpha r_+^{-1}|u_0|_{r,\alpha;s}|\Phi|_{r_+,\alpha;s_+} \ll C_0M\varepsilon_0 r_+^{l-\alpha},\label{u0+est2}\\
& & |u_1(\phi+\Phi)-u_1(\phi)|_{r_+,\alpha;s_+}\leq C_\alpha r_+^{-1}|u_1|_{r,\alpha;s}|\Phi|_{r_+,\alpha;s_+}\ll C_0M\varepsilon_0 r_+^{l-(\alpha+1)(\iota+2)-1},\label{u1+est2}\\
& & |w(\phi+\Phi)-w(\phi)|_{r_+,\alpha;s_+}\leq C_\alpha r_+^{-1}|w|_{r,\alpha;s}|\Phi|_{r_+,\alpha;s_+}\ll C_0M\varepsilon_0 r_+^{l-(\alpha+1)(\iota+2)}.\label{w+est2}
\end{eqnarray}
Lemma~\ref{lemmaA5}, Condition (v.2) and the definition of $K_+$ imply
\begin{equation}\label{u01+est3}
 |({\rm Id}-\Gamma_{K_+}) u_0|_{r_+,\alpha;s_+}\leq r_+ C_0M\varepsilon_0 r_+^{l-\alpha},\quad |({\rm Id}-\Gamma_{K_+}) u_1|_{r_+,\alpha;s_+}\ll C_0M\varepsilon_0 r_+^{l-(\alpha+1)(\iota+2)-1},
 \end{equation}
\begin{equation}\label{w+est3}
 |({\rm Id}-\Gamma_{K_+}) w|_{r_+,\alpha;s_+}\ll C_0M\varepsilon_0 r_+^{l-(\alpha+1)(\iota+2)},
 \end{equation}
 and \eqref{dv1est}, \eqref{dphiest}, Condition (v.2) and \eqref{A+est},
\begin{eqnarray}
& & |P_1^{-1}v_1P_1\widetilde{A}|_{r_+,\alpha;s_+}\leq C_\alpha |P_1^{-1}v_1P_1|_{r_+,\alpha;s_+}
|\widetilde{A}|_{\alpha;s_+}\ll C_0M\varepsilon_0 r_+^{l-(\alpha+1)(\iota+2)-1},\label{u1+est4}\\
& & |P_2^{-1}\partial_\phi\Phi P_2 \tilde{\omega}|_{r_+,\alpha;s_+}\leq C_\alpha |P_2^{-1}\partial_\phi\Phi P_2|_{r_+,\alpha;s_+}|w|_{r,\alpha;s}\ll C_0M\varepsilon_0 r_+^{l-(\alpha+1)(\iota+2)}. \label{w+est4}
\end{eqnarray}
Combining the estimates \eqref{w+est1}, \eqref{w+est2}, \eqref{w+est3}, \eqref{w+est4}, \eqref{G+est} and \eqref{ddest1} for $w^+$, we have
\begin{equation}\label{w+est}
 |w^+|_{r_+,\alpha;s_+}\leq C_0M\varepsilon_0 r_+^{l-(\alpha+1)(\iota+2)}
\end{equation}
by \eqref{w+exp} and choosing small $\varepsilon_0^*$. By the estimates \eqref{dv0est}, \eqref{dv1est}, \eqref{Qest}, \eqref{ddest1}, \eqref{w+est} and Condition (v.2), we also have
\begin{eqnarray}
& &  |P_1^{-1}\partial_\phi v_0 P_2( \tilde{\omega}+w^+)|_{r_+,\alpha;s_+}\leq C_\alpha |P_1^{-1}\partial_\phi v_0 P_2|_{r_+,\alpha;s_+}(|w|_{r,\alpha;s} + |w^+|_{r_+,\alpha;s_+}) \ll C_0M\varepsilon_0 r_+^{l-\alpha},\qquad \label{u0+est4}\\
& & |P_1^{-1}\partial_\phi v_1 P_2( \tilde{\omega}+w^+)|_{r_+,\alpha;s_+} \ll C_0M\varepsilon_0 r_+^{l-(\alpha+1)(\iota+2)-1} \label{u1+est6}
\end{eqnarray}
and
\begin{equation}\label{u1+est5}
|P_1^{-1}\partial_\phi v_0(E_{n_2}+\partial_\phi\Phi)^{-1} P_2Q |_{r_+,\alpha;s_+}\leq C_\alpha |P_1^{-1}\partial_\phi v_0 P_2|_{r_+,\alpha;s_+} |Q|_{r_+,\alpha;s_+} \ll C_0M\varepsilon_0 r_+^{l-(\alpha+1)(\iota+2)-1}.
\end{equation}
From Lemma~\ref{lemmaA6} (ii), Condition (v.2), \eqref{dv0est} and \eqref{dv1est}, it follows
\begin{eqnarray}
& &  |u_1(\phi+\Phi)v_0|_{r_+,\alpha;s_+}\leq C_\alpha |u_1(\phi+\Phi)|_{r_+,\alpha;s_+}|v_0|_{r_+,\alpha;s_+}\ll C_0M\varepsilon_0 r_+^{l-\alpha},\label{u0+est5}\\
& &  |u_1(\phi+\Phi)v_1 |_{r_+,\alpha;s_+}\leq C_\alpha |u_1(\phi+\Phi)|_{r_+,\alpha;s_+} |v_1|_{r_+,\alpha;s_+} \ll C_0M\varepsilon_0 r_+^{l-(\alpha+1)(\iota+2)-1}.\label{u1+est7}
\end{eqnarray}
On account of the estimates \eqref{u0+est1}, \eqref{u1+est1}, \eqref{u0+est2}, \eqref{u1+est2}, \eqref{u01+est3}, \eqref{u0+est4}-\eqref{u1+est7} for $u_0^+$ and $u_1^+$, \eqref{ddest1} and \eqref{G+est}-\eqref{dG+est}, and the expressions \eqref{u0+exp} and \eqref{u1+exp}, we can choose $\varepsilon_0^*$ so small that our estimates yield
$$  |u_0^+|_{r_+,\alpha;s_+}\leq 4 C_0M\varepsilon_0 r_+^{l-\alpha},\qquad
|u_1^+|_{r_+,\alpha;s_+}\leq C_0M\varepsilon_0 r_+^{l-(\alpha+1)(\iota+2)-1}.$$
To turn to the estimates of $H_1^+$ and $H_2^+$, by \eqref{H2+exp} and \eqref{H1+exp} we have
\begin{eqnarray}
 H_1^+-H_1& = & P_1^{-1}(E_{n_1}+v_1)^{-1}P_1\left[H_1\circ T^+(\rho,\phi)-H_1(\rho,\phi)-\partial_I H_1(v_0,\phi+\Phi)(E_{n_1}+v_1)\rho\right.\nonumber\\
& & \,\left.-H_1(v_0,\phi+\Phi)+\mathfrak{G}_1(\rho,\phi)-\mathfrak{G}_1(0,\phi)-\partial_\rho\mathfrak{G}_1(0,\phi)\rho\right.\nonumber\\
& & \left.+P_1^{-1}\partial_\phi v_0(E_{n_2}+\partial_\phi\Phi)^{-1} P_2Q\rho -P_1^{-1}(\partial_\phi v_0+\partial_\phi v_1\rho) P_2H_2^+-P_1^{-1}v_1P_1H_1\right] \label{H1+est}
\end{eqnarray}
and
\begin{eqnarray}
H_2^+-H_2& = & P_2^{-1}(E_{n_2}+\partial_\phi\Phi)^{-1}P_2\left[H_2\circ T^+(\rho,\phi)-H_2(\rho,\phi)-H_2(v_0,\phi+\Phi)\right.\nonumber \\
& & \,\left. +\mathfrak{G}_2(\rho,\phi)-\mathfrak{G}_2(0,\phi)  -P_2^{-1}\partial_\phi\Phi P_2H_2(\rho,\phi)   \right].\label{H2+est}
\end{eqnarray}
After a short calculation, we find
$$|H_1\circ T^+-H_1|_{\mathcal{V}_+,\alpha;s_+} \leq C_\alpha r_+^{-1} |H_1|_{\mathcal{V},\alpha;s}| T^+-{\rm Id}|_{\mathcal{V}_+,\alpha;s_+}\ll \chi_+C_0M\varepsilon_0$$
and
$$|H_2\circ T^+-H_2|_{\mathcal{V}_+,\alpha;s_+} \ll \chi_+C_0M\varepsilon_0$$
by Lemma~\ref{lemmaA6} (i), the Cauchy inequality, \eqref{vTest} and \eqref{Hiest}. It implies
$$|\mathfrak{G}_1(\rho,\phi)-\mathfrak{G}_1(0,\phi)-\partial_\rho\mathfrak{G}_1(0,\phi)\rho|_{\mathcal{V}_+,\alpha;s_+} \ll \chi_+C_0M\varepsilon_0$$
and
$$\quad |\mathfrak{G}_2(\rho,\phi)-\mathfrak{G}_2(0,\phi)|_{\mathcal{V}_+,\alpha;s_+} \ll \chi_+C_0M\varepsilon_0$$
by \eqref{G+est},
$$|P_1^{-1}v_1P_1H_1|_{\mathcal{V}_+,\alpha;s_+} \leq C_\alpha |P_1^{-1}v_1P_1|_{r_+,\alpha;s_+} |H_1|_{\mathcal{V}_+,\alpha;s_+} \ll \chi_+C_0M\varepsilon_0$$
and
$$\quad
|P_2^{-1}\partial_\phi\Phi P_2H_2|_{\mathcal{V}_+,\alpha;s_+} \ll \chi_+C_0M\varepsilon_0$$
by \eqref{dv1est}, \eqref{dphiest} and \eqref{Hiest}, and
$$|P_1^{-1}\partial_\phi v_0(E_{n_2}+\partial_\phi\Phi)^{-1} P_2Q\rho |_{\mathcal{V}_+,\alpha;s_+} \leq C_\alpha |P_1^{-1}\partial_\phi v_0P_2|_{r_+,\alpha;s_+} |Q|_{r_+,\alpha;s_+} \ll \chi_+C_0M\varepsilon_0$$
by \eqref{ddest1}, \eqref{dv0est} and \eqref{Qest}. The above estimates, \eqref{ddest1}, \eqref{w+est1} and the expression \eqref{H2+est} yield
$$|H_2^+-H_2|_{\mathcal{V}_+,\alpha;s_+} \leq \chi_+C_0M\varepsilon_0,\qquad |H_2^+|_{\mathcal{V}_+,\alpha;s_+} \leq (X_++r_0^{l-1})C_0M\varepsilon_0 <\widetilde{C}_0C_0M\varepsilon_0, $$
which, together with \eqref{dv0est} and \eqref{dv1est}, implies
$$|P_1^{-1}(\partial_\phi v_0+\partial_\phi v_1\rho) P_2H_2^+|_{\mathcal{V}_+,\alpha;s_+} \ll \chi_+C_0M\varepsilon_0.$$
Thus, the above estimates, \eqref{ddest1}, \eqref{u0+est1}-\eqref{u1+est1} and the expression \eqref{H1+est} also yield
$$|H_1^+-H_1|_{\mathcal{V}_+,\alpha;s_+} \leq \chi_+C_0M\varepsilon_0,\qquad |H_1^+|_{\mathcal{V}_+,\alpha;s_+} \leq (X_++2r_0^{l-2})C_0M\varepsilon_0<\widetilde{C}_0C_0M\varepsilon_0. $$
Obviously,
$$H_1^+=\varepsilon_0O_{\mathcal{V}_+,\alpha;s_+}(\rho^2),\qquad H_2^+=\varepsilon_0O_{\mathcal{V}_+,\alpha;s_+}(\rho).$$
This completes the proof of the Iteration Lemma. \hskip 0.4in $\blacksquare$

\vskip 0.2in
\renewcommand{\theequation}{\thesection.\arabic{equation}}
\section*{5  Proof of Theorems 2 and 3}
\setcounter{section}{5}\setcounter{equation}{0}

By \eqref{frees}, (H1), (H2) and the Whitney extension theorem (see Lemma~\ref{lemmaA.1}), still using $\omega^\nu$ and $\Lambda^\nu$ to denote their extensions, we have
$$||\Lambda^\nu||_{\alpha;\Pi}\leq 2c_1\varepsilon^{q_3}\leq 2c_1\varepsilon^{q_5},$$
and also
$$\max_{0\leq \mu \leq \alpha}\parallel D^{\mu}\langle b,\omega^\nu(\xi,\varepsilon)\rangle\parallel \geq \frac{c_2}{2}\varepsilon^{q_5} \qquad {\rm in \, Case}\,\, n_{22}=0,$$
\begin{equation}\label{nondeg2}
\max_{1\leq \mu \leq \alpha}\parallel D^{\mu}\langle b,\omega^\nu(\xi,\varepsilon)\rangle\parallel \geq \frac{c_2}{2}\varepsilon^{q_5} \qquad {\rm in \, Case}\,\, n_{22}\neq 0
\end{equation}
by \eqref{Dn11} and \eqref{Dn22}, for sufficiently small $\varepsilon$, all $\xi\in \Pi, b\in \mathcal{S}_{n_2,1}, \nu=0,1,2,\cdots. $

Set
$$f_{km}^\nu(\xi)=\langle k,\omega^\nu(\xi,\varepsilon)\rangle+\langle m,{\rm Im}\Lambda^\nu(\xi,\varepsilon)\rangle$$
for  $0\neq k\in \mathbb{Z}^{n_2}, K_\nu<|k|_2\leq K_{\nu+1};  m={\rm col}(m_1,\cdots,m_{n_1})\in \mathbb{Z}^{n_1}, |m|_1\leq 2$ and $m_1+\cdots+m_{n_1}=0$ or $-1$, $\nu=0,1,2,\cdots $. Here ${\rm Im}\Lambda^\nu$ is the imaginary part of $\Lambda^\nu$. Then
$$\mathfrak{R}_{km}^\nu(\gamma)\subset \left\{\xi\in \Pi_0:\, \mid f_{km}^{\nu-1}(\xi)\mid<\gamma \varepsilon^{q_5}|k|_2^{-\iota}\right\},\quad \nu=1,2,\cdots,$$
where $\Pi_0$ is the closed subset of $\Pi$ defined in Theorem 1.

\begin{lemma}\label{lemma5.1}
If $|k|_2\geq \frac{16}{c_2}c_1|m|_1n_3^{\alpha/2}$, $0\neq k\in \mathbb{Z}^{n_2}$. Then
\begin{equation}\label{measest1}
{\rm meas}\mathfrak{R}_{km}^\nu(\gamma)\leq c_5({\rm diam}\Pi_0)^{n_3-1}(\gamma |k|_2^{-\iota-1})^{\frac{1}{\alpha}}
\end{equation}
for some positive constant $c_5$, where ${\rm diam}\Pi_0$ represents the diameter of $\Pi_0$.
\end{lemma}

{\bf Proof} We only give the proof for the case $n_{22}\neq 0$, the proof of the case $n_{22}=0$ is analogous and is omitted.

Due to the continuity of the derivatives and the compactness of $\Pi_0$ and $\mathcal{S}_{n_2,1}$, the non-degenerate condition \eqref{nondeg2} implies that there exist finite covers $\{\Pi^i\}_{i=1}^{i_0}$ and $\{\mathcal{S}^j\}_{j=1}^{j_0}$ of $\Pi_0$ and $\mathcal{S}_{n_2,1}$, respectively, and $\mu_{ij}: 1\leq \mu_{ij}\leq \alpha, i=1,\cdots, i_0; j=1,\cdots, j_0$, $\Pi^i$ is chosen to be convex, such that
$$\parallel D^{\mu_{ij}}\langle b,\omega^\nu(\xi,\varepsilon)\rangle\parallel \geq \frac{c_2}{4}\varepsilon^{q_5} \qquad {\rm for \, all}\,\,\xi\in\Pi^i, b\in\mathcal{S}^j.$$
Hence, for $ 0\neq k\in \mathbb{Z}^{n_2}, \frac{k}{|k|_2}\in \mathcal{S}^{j_k}$, we have
\begin{equation}\label{nondeg3}
\parallel D^{\mu_{ij_k}}f_{km}^\nu(\xi)\parallel\geq |k|_2\parallel D^{\mu_{ij_k}}\langle \frac{k}{|k|_2},\omega^\nu(\xi,\varepsilon)\rangle\parallel -2c_1|m|_1 \varepsilon^{q_5} n_3^{\frac{\mu_{ij_k}}{2}}\geq \frac{c_2}{8}\varepsilon^{q_5}|k|_2
\end{equation}
for all $\xi\in \Pi^i, i=1,\cdots, i_0$, admitted $m$ and $\nu$ if $|k|_2\geq \frac{16}{c_2}c_1|m|_1n_3^{\alpha/2}$.

Now we estimate the measure of $\mathfrak{R}_{km}^\nu(\gamma)\bigcap \Pi^i$. It follows by \eqref{nondeg3} and the definition of the norm (see Lemma~\ref{lemmaA.3}) that there is a vector $a\in\mathcal{S}_{n_3,1}$ such that
\begin{equation}\label{nondeg4}
\mid D^{\mu_{ij_k}}f_{km}^\nu(\xi)a^{\otimes\mu_{ij_k}}\mid \geq \frac{c_2}{8}\varepsilon^{q_5}|k|_2 \qquad {\rm for \, all}\,\,\xi\in\Pi^i.
\end{equation}
Write $\xi=at+\zeta$ with $t\in \mathbb{R}, \zeta\in a^\perp$ and let $f(t)=f_{km}^\nu(at+\zeta)$, $I_{\zeta}=\{t\in \mathbb{R}:\, at+\zeta\in \Pi^i\}$. The inequality \eqref{nondeg4} means
$$\left|\frac{{\rm d}^{\mu_{ij_k}}f(t)}{{\rm d}t^{\mu_{ij_k}}}\right|\geq \frac{c_2}{8}\varepsilon^{q_5}|k|_2 \qquad {\rm for \, all}\,\,t\in I_{\zeta}.$$
By Fubini's theorem and Lemma~\ref{lemmaA.2}, it implies
$${\rm meas}\left(\mathfrak{R}_{km}^\nu(\gamma)\bigcap \Pi^i\right)\leq c_6({\rm diam}\Pi_0)^{n_3-1}(\gamma |k|_2^{-\iota-1})^{\frac{1}{\alpha}},$$
where $c_6=4\max\{1,(4\alpha !/c_2)^{1/\alpha}\}$. Therefore
$${\rm meas}\mathfrak{R}_{km}^\nu(\gamma)\leq \sum_{i=1}^{i_0}{\rm meas}\left(\mathfrak{R}_{km}^\nu(\gamma)\bigcap \Pi^i\right)\leq i_0 c_6({\rm diam}\Pi_0)^{n_3-1}(\gamma |k|_2^{-\iota-1})^{\frac{1}{\alpha}}.$$
The estimate \eqref{measest1} is proved by setting $c_5=i_0c_6$. \hskip 0.4in $\blacksquare$
\vskip 0.2in
Now, let
$$\mathfrak{K}=\{(k,m)\in \mathbb{Z}^{n_2}\times \mathbb{Z}^{n_1}:\, 0<|k|_2<K^*, 1\leq |m|_1\leq 2, m_1+\cdots+m_{n_1}=0 \,{\rm or} \,-1\}.$$
By an analogous proof to Lemma~\ref{lemma5.1}, we also have

\begin{lemma}\label{lemma5.2}
If $(k,m)\in \mathfrak{K}$ and the condition (ii)' in Theorem 3 holds, then there is a constant $c_7>0$ such that
\begin{equation}\label{measest2}
{\rm meas}\mathfrak{R}_{km}^\nu(\gamma)\leq c_7({\rm diam}\Pi_0)^{n_3-1}(\gamma |k|_2^{-\iota})^{\frac{1}{\alpha}}.
\end{equation}
\end{lemma}

\vskip 0.2in
{\bf Remark 5.1} If $q_3>q_5\geq 0$, then without the conditions (ii) and (ii)', we can obtain
$${\rm meas}\mathfrak{R}_{km}^\nu(\gamma)\leq c_5({\rm diam}\Pi_0)^{n_3-1}(2\gamma |k|_2^{-\iota-1})^{\frac{1}{\alpha}}\quad {\rm for\, all}\,\ (k,m)\in \mathfrak{K}.$$

In fact, for sufficiently small $\gamma$ (equivalently, sufficiently small $\varepsilon$), we have
$$\mathfrak{R}_{km}^\nu(\gamma)\subset \left\{\xi\in \Pi_0:\, |\langle k, \omega^\nu(\xi,\varepsilon)\rangle|< 2\gamma \varepsilon^{q_5}|k|_2^{-\iota}\right\}$$
for all $(k,m)\in \mathfrak{K}$. From  Lemma~\ref{lemma5.1} with $|m|_1=0$, it follows
\begin{eqnarray*}
{\rm meas}\mathfrak{R}_{km}^\nu(\gamma)&\leq & {\rm meas} \left\{\xi\in \Pi_0:\, |\langle k, \omega^\nu(\xi,\varepsilon)\rangle|< 2\gamma \varepsilon^{q_5}|k|_2^{-\iota}\right\}\\
&\leq &  c_5({\rm diam}\Pi_0)^{n_3-1}(2\gamma |k|_2^{-\iota-1})^{\frac{1}{\alpha}}.
\end{eqnarray*}

\vskip 0.2in
{\bf Proof of Theorem 3}
By Lemmas ~\ref{lemma5.1} and ~\ref{lemma5.2}, we obtain
\begin{eqnarray*}
{\rm meas}(\Pi_0\setminus \Pi_{\gamma})&\leq & \sum_{\nu=1}^{\infty}\sum_{K_{\nu-1}<|k|_2\leq K_\nu}\left({\rm meas}\mathfrak{R}_{k0}^\nu +\sum_{|m|_1=1} {\rm meas}\mathfrak{R}_{km}^\nu + \sum_{|m|_1=2} {\rm meas}\mathfrak{R}_{km}^\nu\right)\\
& \leq & \sum_{(k,m)\in \mathfrak{K}} {\rm meas}\mathfrak{R}_{km}^\nu + c_5({\rm diam}\Pi_0)^{n_3-1} \gamma^{\frac{1}{\alpha}} \left[ \sum_{0\neq k\in \mathbb{Z}^{n_2}} |k|_2^{-\frac{\iota+1}{\alpha}}\right.\\
& & + \left. n_1 \sum_{|k|_2\geq K^*} |k|_2^{-\frac{\iota+1}{\alpha}} + n_1(n_1-1) \sum_{|k|_2\geq K^*} |k|_2^{-\frac{\iota+1}{\alpha}}\right]\\
& \leq & c_{10}\gamma^{\frac{1}{\alpha}},
\end{eqnarray*}
where $c_{10}$ is a positive constant depending on $n_1, n_2, {\rm diam}\Pi, \omega_i$ and $\Lambda_i (i=1,2)$, and ${\rm meas}(\Pi\setminus \Pi_0)=O(\gamma)$. This proves Theorem 3.
 \hskip 0.4in $\blacksquare$

\vskip 0.2in
{\bf Proof of Theorem 2}  From the proof of Theorem 3, it is easy to see
$${\rm meas}(\bigcup_{\nu=1}^{\infty}\bigcup_{(k,m)\notin \mathfrak{K}} \mathfrak{R}_{km}^\nu(\gamma))\leq c_{10}\gamma^{\frac{1}{\alpha}}\rightarrow 0 \qquad {\rm as}\, \gamma \rightarrow 0.$$
On the other hand, there is a $\nu_0$ such that $K_{\nu_0}\leq K^*$, hence
$$\bigcup_{\nu=1}^{\infty}\bigcup_{(k,m)\in \mathfrak{K}} \mathfrak{R}_{km}^\nu(\gamma)\subset \bigcup_{\nu=1}^{\nu_0}\bigcup_{(k,m)\in \mathfrak{K}} \mathfrak{R}_{km}^\nu(\gamma).$$
By the condition (ii) in Theorem 2, the boundedness of $\Pi$, and \eqref{I2}, we have ${\rm meas}\mathfrak{R}_{km}^\nu(\gamma)\rightarrow 0$ as $\gamma\rightarrow 0$. Since $ \mathfrak{K}$ is finite, we also have
$${\rm meas}(\bigcup_{\nu=1}^{\nu_0}\bigcup_{(k,m)\in \mathfrak{K}} \mathfrak{R}_{km}^\nu(\gamma))\rightarrow 0\qquad  {\rm as}\, \gamma\rightarrow 0. \hskip 0.4in \blacksquare$$

\vskip 0.2in
\renewcommand{\theequation}{\thesection.\arabic{equation}}
\section*{ Appendix}
\setcounter{section}{6}\setcounter{equation}{0}
\vskip 0.1in
\subsection*{A.1. Whitney extension theorem }
\vskip 0.1in
Let $\overline{\Omega}\subset \mathbb{R}^n$ be a closed set, $p$ be a non-negative integer, $p<l\leq p+1$. $C_W^l(\overline{\Omega})$ is the class of all collections $f=\{f^{(k)}\}_{|k|_1\leq p}$ of functions defined on $\overline{\Omega}$ which satisfy, for some finite $M$,
\begin{equation}\label{whitineq}
|f^{(k)}(x)|\leq M,\qquad |f^{(k)}(x)-P_k(x,y)|\leq M |x-y|^{l-|k|_1}
\end{equation}
for all $x,y\in \overline{\Omega}$ and $|k|_1\leq p$, where
$$P_k(x,y)=\sum_{|k+j|_1\leq p}\frac{1}{j!}f^{(k+j)}(y)(x-y)^j$$
is the analogue of the $k$-th Taylor polynomial. $f$ is called $C^l$ Whitney in $\overline{\Omega}$ with Whitney derivatives $D^k f=f^{(k)}$ for $|k|_1\leq p$. Define a norm
$$||f||_{C_W^l(\overline{\Omega})}=\inf M$$
is the smallest $M$ for which both inequalities in \eqref{whitineq} hold. Then $C_W^l(\overline{\Omega})$ with the norm is a Banach space.

The following extension theorem indicates that a Whitney differentiable function has an extension to $\mathbb{R}^n$ which is differentiable in the standard sense.

\begin{lemma}\label{lemmaA.1}
(Whitney extension theorem, \cite{Whit34, Stein70, Pos82}) Let  $\overline{\Omega}$ be a closed set in $\mathbb{R}^n$, $p\in \mathbb{Z}_+$ and $p<l\leq p+1$. Then there exists a linear extension operator
$$\mathfrak{E}: C_W^l(\overline{\Omega})\rightarrow C^l(\mathbb{R}^n), \qquad f=\{f^{(k)}\}_{|k|_1\leq p}\rightarrow F=\mathfrak{E}f$$
such that
$$D^kF\mid_{\overline{\Omega}}=f^{(k)}, \qquad |k|_1\leq p$$
and
$$||F||_{l;\mathbb{R}^n}\leq C ||f||_{C_W^l(\overline{\Omega})},$$
where the constant $C$ depends only on $l$ and the dimension $n$, but not on $\overline{\Omega}$. Moreover, if $\overline{\Omega}=\overline{\Omega}_1\times \mathbb{T}^{n_2}\subset \mathbb{R}^{n_1}\times \mathbb{T}^{n_2}$, then the extension can be chosen to be defined on $\mathbb{R}^{n_1}\times \mathbb{T}^{n_2}$, so that the periodicity is preserved.
\end{lemma}

\vskip 0.2in
\subsection*{A.2. Measure estimate lemmas}
\vskip 0.1in
\begin{lemma}\label{lemmaA.2}
\cite{Russ01} Let $f:[a,b]\rightarrow \mathbb{R}$ with $a<b$ be an $\alpha$-times continuously differentiable function satisfying
$$\left| \frac{{\rm d}^\alpha f(x)}{{\rm d}x^\alpha}\right| \geq c, \qquad x\in [a,b]$$
for some $\alpha \in \mathbb{N}$ and a constant $c>0$. Then we have the measure estimate
$${\rm meas}\{x\in [a,b]:\, |f(x)|\leq \varepsilon\}\leq 4\left(\alpha !\frac{\varepsilon}{2c}\right)^{\frac{1}{\alpha}} \quad {\rm for\, all}\,\,\varepsilon>0.$$
\end{lemma}

\begin{lemma}\label{lemmaA.3}
Let $\Pi\subset \mathbb{R}^p$ be a bounded closed set, $f_j: \Pi\rightarrow \mathbb{R}$ be of $C^{\alpha}$ on $\Pi$ with a positive integer $\alpha, j=1,\cdots,q$. Denote $f(\xi)={\rm col}(f_1(\xi),\cdots,f_q(\xi))$. Assume for $\xi\in \Pi$,
\begin{equation}\label{rankcon}
{\rm rank}\left( f(\xi), \frac{\partial^\beta f(\xi)}{\partial \xi^\beta}:\, 1\leq |\beta|_1\leq \alpha\right)=q \quad {\rm and} \quad 1+p+p^2+\cdots+p^{\alpha}\geq q.
\end{equation}
Then there is a constant $c>0$ such that
$$\max_{0\leq \mu\leq \alpha}||D^\mu\langle b, f(\xi)\rangle||\geq c \qquad {\rm for\, all}\,\,b\in \mathcal{S}_{q,1},\, \xi\in\Pi.$$
Here $D$ represents the differential operator with respect to the variable $\xi$,
$$\mathcal{S}_{q,1}=\{b\in \mathbb{R}^q:\, |b|_2=1\}, \quad ||D^\mu\langle b, f(\xi)\rangle||=\max_{a\in\mathcal{S}_{p,1}} |D^\mu\langle b, f(\xi)\rangle a^{\otimes\mu}|,$$
$a^{\otimes\mu}=(a_1,a_2,\cdots,a_\mu)$ with $a_i=a, i=1,2,\cdots,\mu$.
\end{lemma}

{\bf Remark A.1} Here, by the Whitney extension theorem  we assume the continuous differentiability of a function $f$ with respect to the parameter variable $\xi$ on a closed set $\Pi$ means that $f$ is continuously differentiable in some neighbourhood of $\Pi$.

{\bf Proof} Suppose such a constant $c$ does not exist. Then for any positive integer $n$, we can find $\xi_n\in\Pi$ and $b_n\in \mathcal{S}_{q,1}$ satisfying
$$\max_{0\leq \mu\leq \alpha}||D^\mu\langle b_n, f(\xi_n)\rangle||<\frac{1}{n},\qquad n=1,2,\cdots.$$
Based on the compactness of $\Pi$ and $\mathcal{S}_{q,1}$ there are convergent subsequences of $\{b_n\}$ and $\{\xi_n\}$, respectively, still denoting by $\{b_n\}$ and $\{\xi_n\}$, such that $b_n\rightarrow b_0\in \mathcal{S}_{q,1}, \xi_n\rightarrow \xi_0\in \Pi$ as $n\rightarrow \infty$. Thus, the continuity of the derivatives implies
$$||D^\mu\langle b_0, f(\xi_0)\rangle||=0 \qquad {\rm for\, all}\,\,0\leq \mu\leq \alpha.$$
Noting that
$$|D^\mu\langle b_0, f(\xi_0)\rangle (a_1,\cdots,a_\mu)|\leq \frac{\mu^\mu}{\mu !}||D^\mu\langle b_0, f(\xi_0)\rangle||$$
for all $a_i\in \mathcal{S}_{p,1},i=1,\cdots,\mu$, we have
$$b_0^T\left( f(\xi_0), \frac{\partial^\beta f(\xi_0)}{\partial \xi^\beta}:\, 1\leq |\beta|_1\leq \alpha\right)=0,$$
which implies
$${\rm rank}\left( f(\xi_0), \frac{\partial^\beta f(\xi_0)}{\partial \xi^\beta}:\, 1\leq |\beta|_1\leq \alpha\right)<q$$
being in contradiction with the condition \eqref{rankcon}. The lemma is proved. \hskip 0.4in $\blacksquare$

{\bf Remark A.2} From the proof of Lemma~\ref{lemmaA.3}, it is easy to see that if the condition \eqref{rankcon} is replaced by
$${\rm rank}\left( \frac{\partial^\beta f(\xi)}{\partial \xi^\beta}:\, 1\leq |\beta|_1\leq \alpha\right)=q \qquad {\rm and} \quad p+p^2+\cdots+p^\alpha\geq q,$$
then we also have
$$\max_{1\leq \mu\leq \alpha}||D^\mu\langle b, f(\xi)\rangle||\geq c \qquad {\rm for\, all}\,\,b\in \mathcal{S}_{q,1},\, \xi\in\Pi.$$

\vskip 0.2in
\subsection*{A.3. Properties of analytic smoothing operator}
\vskip 0.1in
Let $l>0, m\in \mathbb{N}$ and $C^l(\mathbb{R}^m)$ be the H\"{o}lder space defined in Definition~\ref{def1} without parameter variables, $u_0\in C_0^\infty(\mathbb{R})$ be an even function, vanishing outside the interval $[-1,1]$ and identically equal to 1 in a neighbourhood of 0, $u(x)=u_0(|x|_2^2)$ for $x\in \mathbb{R}^m$ and
$$\tilde{u}(z)=\int_{\mathbb{R}^m} u(x)e^{\sqrt{-1}\langle z,x\rangle}dx \qquad {\rm for} \,\, z\in \mathbb{C}^m,$$
$$f_r(x)\coloneqq (\mathcal{S}_r f)(x)\coloneqq r^{-m}\int_{\mathbb{R}^m} \tilde{u}((x-y)/r)f(y)dy$$
for $x\in \mathbb{C}^m$ and $r\in (0,1]$.

\begin{lemma}\label{lemmaA4}
The following assertions are valid

{\rm (a)}\, $\int_{\mathbb{R}^m} \tilde{u}(x) dx=u(0)=1;$

{\rm (b)}\, $(\mathcal{S}_r f)(x)=\int_{\mathbb{R}^m} \tilde{u}(y)f(x-ry)dy$ \, for  $x\in \mathbb{R}^m$;

{\rm (c)}\, $\int_{\mathbb{R}^m} x^k \tilde{u}(x) dx=0$ \,  for $0\neq k\in \mathbb{Z}_+^m$;

{\rm (d)}\, for any $p\in \mathbb{N}$, there is a constant $C_p>0$ such that
$$\left|D^k\tilde{u}(z)\right|\leq \frac{C_p}{(1+|z|_2)^p} e^{|{\rm Im}z|_2} \qquad {\rm for \, all}\,\, |k|_1\leq p, k\in \mathbb{Z}_+^m,$$
\hskip 0.4in where $D^k=D_1^{k_1}\circ D_2^{k_2}\circ \cdots \circ D_{m}^{k_m}$, and $D_j^{k_j}=\frac{\partial^{k_j}}{\partial x_j^{k_j}}$;

{\rm (e)}\, if $P$ is a polynomial, then $(\mathcal{S}_r P)(x)=P(x)$;

{\rm (f)}\, there exists a constant $C_l>0$ such that
$$\left|D^kf_r(x)-\sum_{|\beta|_1\leq l-|k|_1} D^{k+\beta}f({\rm Re} x)\frac{(\sqrt{-1}{\rm Im} x)^\beta}{\beta!}\right|\leq C_l r^{l-|k|_1}||f||_{l;\mathbb{R}^m},\quad |{\rm Im } x|\leq r \leq 1$$
for all $k\in \mathbb{Z}_+^m$ with $|k|_1\leq l$. In particular, for $x\in \mathbb{R}^m$ and $p\in \mathbb{Z}_+$,
$$||f_r-f||_{p;\mathbb{R}^m}\leq C_{lp} r^{l-p}||f||_{l;\mathbb{R}^m},\qquad p\leq l$$
for a suitable constant $C_{lp}$ depending on $l,p$ and $m$.
\end{lemma}

{\bf Proof} The definitions of $\tilde{u}$ and $\mathcal{S}_r$ imply (a) and (b), respectively. Noting that the  $\tilde{u}$ is a Schwartz function (see (d)), and the Fourier transformation and differentiation can be exchanged, we have
$$ \int_{\mathbb{R}^m} x^k \tilde{u}(x) dx=\left.(\sqrt{-1})^{|k|_1}D_y^k\int_{\mathbb{R}^m}  \tilde{u}(x)e^{-\sqrt{-1}\langle y,x\rangle} dx\right|_{y=0}= \left.(\sqrt{-1})^{|k|_1}D^ku(y)\right|_{y=0}=0,$$
which verifies (c). See Lemma 9, Proposition 8 and Remark 15 (i) in \cite{Chier03} for (d), (e) and (f), respectively, also see the proof of Lemma 2.1 in Part I of \cite{Zehn75} for (e).
\hskip 0.4in $\blacksquare$

\begin{lemma}\label{lemmaA5}
Let $K$ be a positive integer and $f$ be a bounded and analytic function in the strip $\{x:\,|{\rm Im} x|<r\}$ of $\mathbb{T}^n$, $f(x)=\sum_{k\in \mathbb{Z}^n} \hat{f}(k) e^{\sqrt{-1}\langle k,x\rangle}$. Define the truncation operator $\Gamma_K$ as follows
$$\Gamma_K f=\sum_{|k|_2\leq K}\hat f(k)e^{\sqrt{-1}\langle k,x\rangle}.$$
 If $K>(2\rho)^{-1}$, then we have
 $$|({\rm Id}-\Gamma_K)f|_{r-2\rho} \leq C(n) |f|_r\rho^{-n}e^{-\rho K},\qquad 0<2\rho \leq r,$$
 where $C(n)=6(n!)n^ne^{-n}$.
\end{lemma}

{\bf Proof} Set $\sigma=2\rho$. Based on the fact that the number of all $k$ with $|k|_1=m$ is bounded by $2nm^{n-1}$, we have
\begin{equation}\label{trun1}
|({\rm Id}-\Gamma_K)f|_{r-\sigma} \leq \sum_{|k|_2> K}|f|_r e^{-\sigma |k|_1} \leq \sum_{|k|_1> K}|f|_r e^{-\sigma |k|_1} \leq |f|_r\sum_{m> K}2nm^{n-1} e^{-\sigma m}.
\end{equation}
Here we use Lemma A.1 in \cite{Poschel01}. Since the function $y^{n-1}e^{-\sigma y}$ is monotonically decreasing in the interval $[\frac{n-1}{\sigma},+\infty)$ and $K>\sigma^{-1}$, therefore,
\begin{eqnarray*}
\sum_{m> K}m^{n-1} e^{-\sigma m}<\int_{K}^{+\infty}y^{n-1}e^{-\sigma y}dy & = & \left(\frac{1}{\sigma}K^{n-1}+\frac{n-1}{\sigma^2}K^{n-2}+\cdots +\frac{(n-1)!}{\sigma^n}\right)e^{-\sigma K}\\
& < & 3(n-1)!K^ne^{-\sigma K}.
\end{eqnarray*}
Hence, by \eqref{trun1} we obtain
\begin{equation}\label{trun2}
|({\rm Id}-\Gamma_K)f|_{r-\sigma} \leq 6 n! |f|_rK^n e^{-\sigma K}.
\end{equation}
Noting that the maximum of the function $y^ne^{-y}$ on the interval $(0,+\infty)$ is $n^ne^{-n}$, \eqref{trun2} implies
$$|({\rm Id}-\Gamma_K)f|_{r-2\rho} \leq 6(n!)n^ne^{-n} |f|_r\rho^{-n}e^{-\rho K}. \hskip 0.4in \blacksquare$$
\vskip 0.2in

Let $\Omega_1$ and $\Omega_2$ be domains in $\mathbb{C}^n$, $\Pi$ be an open set in $\mathbb{R}^m$, $f(x,\xi)$ and $g(x,\xi)$ be analytic in $x\in (\Omega_1+r)$ and in $x\in \Omega_2$ respectively, and continuously differential up to order $\alpha$ in $\xi\in \Pi$, $g:\, \Omega_2\times \Pi\rightarrow \Omega_1$, where $r>0,\, \Omega_1+r=\{x\in \mathbb{C}^n:\, {\rm dist}(x,\Omega_1)<r\}$.

We introduce the notation for $1\leq \mu \leq \alpha$,
$$|Df|_\mu \coloneqq \max_{1\leq |\beta|_1\leq \mu}\left|\partial_\xi^\beta Df\right|_{\Omega_1+r, 0; \Pi}, \qquad |g|_\mu \coloneqq \max_{1\leq |\beta|_1\leq \mu}\left|\partial_\xi^\beta g\right|_{\Omega_2, 0; \Pi},$$
where $Df$ represents the differential operator with respect to the coordinate variable $x$. Using the Chain Rule on differentiation of a composition of mappings and Cauchy inequality, we easily prove the following lemma.

\begin{lemma}\label{lemmaA6}
Let $\beta \in \mathbb{Z}_+^m$ and $|\beta|_1=\mu, 1\leq \mu \leq \alpha$. Then

{\rm (i)}\, $\left|\partial_\xi^\beta (f\circ g-f)\right|_{\Omega_2, 0; \Pi}\leq
\left\{\begin{array}{l}\left|(\partial_\xi^\beta f)\circ g-\partial_\xi ^\beta f\right|_{\Omega_2, 0; \Pi}+| Df|_{\Omega_1+r, 0; \Pi}|g|_1 \quad  {\rm for} \quad \mu=1\\
\left|(\partial_\xi^\beta f)\circ g-\partial_\xi ^\beta f\right|_{\Omega_2, 0; \Pi}+| Df|_{\Omega_1+r, 0; \Pi}(|g|_\mu+\frac{(\mu-1)!}{r^{\mu-1}}|g|_1^\mu)\\
\,\, +\sum_{j=1}^{\mu-1}\frac{C_j}{r^{j-1}}|Df|_{\mu-1}|g|_{\mu-1}^j \qquad {\rm for } \quad 2\leq \mu \leq \alpha;
\end{array}\right.$

{\rm (ii)}\,  $\left|\partial_\xi^\beta (f\circ g)\right|_{\Omega_2, 0; \Pi}\leq \left|(\partial_\xi^\beta f)\circ g\right|_{\Omega_2, 0; \Pi}
+\sum_{j=1}^{\mu}\frac{C_j}{r^{j}}|f|_{\mu-1}|g|_{\mu}^j$,\\
where $C_j(j=1,\cdots \mu)$ are nonnegative constants only depending on $\beta$.
\end{lemma}

\vskip 0.2in
\subsection*{A.4. An estimate lemma for small divisors }
\vskip 0.1in
\begin{lemma}\label{lemmaA7}
Assume the frequency vector $\omega=(\omega_1,\cdots, \omega_n)$ satisfies the inequalities
\begin{equation}\label{smdivcon}
|\langle k, \omega\rangle|\geq \frac{\gamma}{|k|_2^\tau} \qquad {\rm and} \qquad |\langle k, \omega\rangle+ \lambda|\geq \frac{\gamma}{|k|_2^\tau}
\end{equation}
for all integer vectors $0\neq k \in \mathbb{Z}^n$ with $|k|_2\leq K \leq \infty$, and some constants $\tau>n-1 \geq 1$, $K>0, \gamma>0$ and $\lambda \in \mathbb{R}$. Then the following inequalities hold
\begin{equation}\label{smdivest0}
\sum_{0\neq |k|_2\leq K} |k|_1^v|\langle k, \omega\rangle|^{-b}e^{-\sigma |k|_1}\leq C \gamma^{-b}\sigma^{-(\tau b+v+1)}
\end{equation}
and
\begin{equation}\label{smdivest}
\sum_{0\neq |k|_2\leq K} |k|_1^v|\langle k, \omega\rangle + \lambda|^{-b}e^{-\sigma |k|_1}\leq C \gamma^{-b}\sigma^{-(\tau b+v+1)}
\end{equation}
with
$$C=15 \tau \sqrt{\tau b +v} 2^{2(n+b)-3} n^{\tau b +v+1}(\tau b -n+1)^{-1}\left(\frac{\tau b +v}{e}\right)^{\tau b+v},$$
where $v\geq 0, b\geq 1$ and $\sigma\in (0,1)$ are constants.
\end{lemma}

{\bf Proof}  The proof is based on the fact that only a few of the denominators $\langle k, \omega\rangle$ and $\langle k, \omega\rangle+ \lambda$ are small, which was used by Siegel \cite{Sieg}, Arnol'd \cite{Arn63} and Moser \cite{Mos66}. For the sake of completeness, we present the proof for our situation and make the involved constants explicit in estimates. We only prove the inequality \eqref{smdivest} with $K=\infty$. The proof of \eqref{smdivest0} and the case $K<\infty$ is analogous and is omitted. Set
$$K(m,j)=\{k=(k_1,\cdots,k_n) \in \mathbb{Z}^n:\, |k|\coloneqq \max_{1\leq i \leq n}|k_i|=m,\, \gamma^{-1}2^j<|\langle k, \omega\rangle+ \lambda|^{-1}\leq \gamma^{-1}2^{j+1}\}$$
and let $K(m,j)^\#$ denote the number of points in $K(m,j)$. Then we have
\begin{equation}\label{Kindex}
K(m,j)^\#\leq (2n)^n(2m)^{n-1}2^{-\frac{1}{\tau}(n-1)(j-1)}.
\end{equation}
In fact, if $k,k'\in K(m,j)$ are different points, then
$$\gamma |k-k'|_2^{-\tau}\leq |\langle k-k', \omega\rangle|\leq |\langle k, \omega\rangle+ \lambda|+|\langle k', \omega\rangle+ \lambda|< \gamma 2^{1-j},$$
which implies
$$|k-k'|\geq n^{-\frac{1}{2}}|k-k'|_2 > n^{-1}2^{\frac{j-1}{\tau}}\coloneqq 2\rho_j.$$
Noting $|k-k'|\leq 2m$ we get $\rho_j\leq m$. If we encircle every point $k\in K(m,j)$ by a cube $\mathfrak{C}_k:\, |x-k|\leq \rho_j$, then these cubes are mutually disjoint. The intersections of these cubes $\mathfrak{C}_k$ with the curved surface $|x|=m$ are disjoint $n-1$ dimensional sets with $n-1$ dimensional volume $\geq \rho_j^{n-1}$. As the $n-1$ dimensional volume of the curved surface $|x|=m$ is $2n(2m)^{n-1}$, we obtain
$$K(m,j)^\#\leq \frac{2n(2m)^{n-1}}{\rho_j^{n-1}},$$
which verifies the inequality \eqref{Kindex}. Thus we have
$$\sum_{K(m,j)}|\langle k, \omega\rangle + \lambda|^{-b}\leq \gamma^{-b}2^{b(j+1)}K(m,j)^\#\leq 2^{2(n+b)-1}n^n m^{n-1}\gamma^{-b}2^{(b-\frac{n-1}{\tau})(j-1)}\coloneqq C_1 2^{(b-\frac{n-1}{\tau})(j-1)}.$$
Let $j^*$ be the greatest occurring $j$ for which $K(m,j)\neq \emptyset$. Then the facts that
$$\gamma^{-1}2^{j^*}<|\langle k, \omega\rangle+ \lambda|^{-1}\leq \gamma^{-1}|k|_2^\tau \leq \gamma^{-1}(nm)^\tau $$
 and
 $$\{ k \in \mathbb{Z}^n:\, |k|=m\}^\#=(2m+1)^n-(2m-1)^n<2n(4m)^{n-1}$$
 imply
 \begin{eqnarray*}
\sum_{|k|=m}|\langle k, \omega\rangle + \lambda|^{-b} & \leq & \sum_{|k|=m,|\langle k, \omega\rangle + \lambda|^{-1}\leq 2\gamma^{-1}}|\langle k, \omega\rangle + \lambda|^{-b}+\sum_{j=1}^{j^*}\sum_{K(m,j)}|\langle k, \omega\rangle + \lambda|^{-b}\\
& \leq & n 2^{2n+b-1}m^{n-1}\gamma^{-b}+C_1\sum_{j=1}^{j^*}2^{(b-\frac{n-1}{\tau})(j-1)}\\
& \leq & \tau 2^{2(n+b)-1}n^{\tau b +1}(\tau b -n+1)^{-1}\gamma^{-b}m^{\tau b}\coloneqq C_2m^{\tau b}.
\end{eqnarray*}
Therefore,
\begin{eqnarray}
& & \sum_{0\neq k\in \mathbb{Z}^n} |k|_1^v|\langle k, \omega\rangle + \lambda|^{-b}e^{-\sigma |k|_1} \leq \sum_{m=1}^\infty \sum_{|k|=m} (nm)^v|\langle k, \omega\rangle + \lambda|^{-b}e^{-\sigma m}\nonumber \\
& & \qquad \leq C_2n^v\sum_{m=1}^\infty m^{\tau b+v}e^{-\sigma m}. \label{sum1}
\end{eqnarray}
Noting that the function $g(x)=x^{\tau b+v}e^{-\sigma x}$ on the interval $[1,\infty)$ gets its maximum at $x_0=\frac{\tau b +v}{\sigma}$, moreover is strictly increasing and decreasing on $[1,x_0)$ and $(x_0,\infty)$, respectively. Denote the integer part of $\frac{\tau b +v}{\sigma}$ by $m_0$. Then $m_0\geq 1$ and
\begin{eqnarray}
\sum_{m=1}^\infty m^{\tau b+v}e^{-\sigma m} & \leq & \int_1^{m_0}x^{\tau b+v}e^{-\sigma x}dx+ g(\frac{\tau b +v}{\sigma}) +\int_{m_0}^\infty x^{\tau b+v}e^{-\sigma x}dx \nonumber\\
& \leq & g(\frac{\tau b +v}{\sigma}) + \sigma^{-(\tau b +v+1)}\int_{0}^\infty y^{\tau b+v}e^{-y}dy \nonumber\\
& = & \left(\frac{\tau b+v}{e\sigma}\right)^{\tau b +v}+\sigma^{-(\tau b +v+1)}\Gamma(\tau b +v+1).\label{sum2}
\end{eqnarray}
By the Stirling formula of the gamma function, we have
\begin{equation}\label{gammaest}
\Gamma(\tau b +v+1)<\frac{11}{4}\sqrt{\tau b +v}\left(\frac{\tau b+v}{e}\right)^{\tau b +v}.
\end{equation}
Combining \eqref{sum1}-\eqref{gammaest}, we obtain the estimate \eqref{smdivest}. The proof of the lemma is complete. \hskip 0.2in $\blacksquare$
\vskip 0.2in

{\bf Remark A.3} From the proof of Lemma~\ref{lemmaA7} it is easy to see that if the norm $|k|_2$ in the condition \eqref{smdivcon} is replaced by the norm $|k|_1$, then the estimates \eqref{smdivest0} and \eqref{smdivest} are still valid.

\end{document}